\documentclass[english]{svmult}

\usepackage{babel}
\usepackage{amstext}
\usepackage{amsmath}
\usepackage{amsfonts}
\usepackage{latexsym}
\usepackage{ifthen}
\usepackage{xypic}
\xyoption{all}

\usepackage{url}
\usepackage{makeidx}         
\usepackage{graphicx}        
\usepackage{multicol}        
\usepackage[bottom]{footmisc}

\makeindex             

\usepackage{eepic,epic}
\newcommand{\N}{{\mathbb N}}

\newcommand{\Z}{{\mathbb Z}}
\newcommand{\R}{{\mathbb R}}
\newcommand{\C}{{\mathbb C}}
\newcommand{\Q}{{\mathbb Q}}
\renewcommand{\P}{{\mathbb P}}

\newcommand{\K}{\C}   

\newcommand{\kh}{{\mathcal H}}

\newcommand{\kj}{{\mathcal J}}

\newcommand{\km}{{\mathcal M}}

\newcommand{\ko}{{\mathcal O}}
\newcommand{\kp}{{\mathcal P}}

\newcommand{\ks}{{\mathcal S}}

\newcommand{\fm}{\mathfrak{m}}

\DeclareMathAlphabet{\mathsc}{U}{rsfs}{m}{n}

\newcommand{\sC}{\mathsc{C}}

\newcommand{\kt}{\mathcal{T}}
\newcommand{\sU}{\mathsc{U}}

\newtheorem{lemma and definition}[theorem]{Lemma and Definition}

\newcommand{\balpha}{\boldsymbol{\alpha}}

\newcommand{\bgamma}{\boldsymbol{\gamma}}

\newcommand{\bH}{\boldsymbol{H}}

\newcommand{\bm}{\boldsymbol{m}}
\newcommand{\bn}{\boldsymbol{n}}

\newcommand{\bone}{\boldsymbol{1}}

\newcommand{\bw}{\boldsymbol{w}}
\newcommand{\bx}{\boldsymbol{x}}

\newcommand{\del}{{\delta}}

\newcommand{\gam}{{\gamma}}

\newcommand{\gea}{\gamma^{\text{\it ea}}}
\newcommand{\ges}{\gamma^{\text{\it es}}}
\newcommand{\Hilb}{\kh ilb\!\;}

\newcommand{\Iea}{I^{\text{\it ea}}}

\newcommand{\Ies}{I^{\text{\it es}}}
\newcommand{\Iesf}{I^{\text{\it es}}_{\text{\it fix}}}

\newcommand{\Ieaf}{I^{\text{\it ea}}_{\text{\it fix}}}
\newcommand{\Isqh}{I^{\text{\it sqh}}}

\newcommand{\lra}{\longrightarrow}

\newcommand{\Sig}{{\Sigma}}
\newcommand{\tci}{\tau_{\text{\it ci}}}

\newcommand{\tes}{\tau^{\text{\it es}}}
\newcommand{\teaci}{\tau_{\text{\it ci}}}
\newcommand{\tesci}{\tau^{\text{\it es}}_{\text{\it ci}}}

\newcommand{\VC}{V_{|C|}(S_1,\dots,S_r)}
\newcommand{\VCf}{V_{|C|,\text{\it fix}}}
\newcommand{\Vd}{V_d^{\text{\it irr}} (S_1,\dots,S_r)}
\newcommand{\Vdi}{V_d^{\text{\it irr}}}
\newcommand{\Vdnk}{\mbox{$V_d^{\text{\it irr}}(n\!\!\:\cdot\!\!\:A_1,
k\!\!\:\cdot \!\!\:A_2)$}}
\newcommand{\VD}{\mbox{$V_{|D|} (S_1,\dots,S_r)$}}
\newcommand{\VDi}{\mbox{$V_{|D|}^{\text{\it irr}} (S_1,\dots,S_r)$}}
\newcommand{\Vhf}{V_{h,\text{\it fix}}}
\newcommand{\Vhi}{\mbox{$V^{\text{\it irr}}_h$}}

\newcommand{\Zast}{Z^{\text{\it a}}_{\text{\it st}}}

\newcommand{\Zea}{Z^{\text{\it ea}}}
\newcommand{\Zes}{Z^{\text{\it es}}}
\newcommand{\Zeaf}{Z^{\text{\it ea}}_{\text{\it fix}}}
\newcommand{\Zesf}{Z^{\text{\it es}}_{\text{\it fix}}}
\newcommand{\Zf}{Z_{\text{\it fix}}}

\newcommand{\Zsqh}{Z^{\text{\it sqh}}}
\newcommand{\Zsst}{Z^{\text{\it s}}_{\text{\it st}}}

\DeclareMathOperator{\const}{const}

\DeclareMathOperator{\mt}{mt}

\DeclareMathOperator{\pr}{pr}

\DeclareMathOperator{\Sing}{Sing}

\DeclareMathOperator{\Sym}{Sym}

\DeclareMathOperator{\Tor}{Tor}
\DeclareMathOperator{\tHilb}{Hilb\!\:}

\DeclareMathOperator{\wdeg}{\bw-deg}

\begin{document}

\title*{Equisingular Families of Projective Curves}
\author{Gert-Martin Greuel \and Christoph Lossen \and Eugenii Shustin}
\authorrunning{G.-M.\ Greuel, C.\ Lossen, and E. Shustin}
\institute{TU Kaiserslautern,
Fachbereich Mathematik,
Erwin-Schr\"odinger-Stra\ss e,
67663 Kaiserslautern, Germany, \texttt{greuel@mathematik.uni-kl.de} \and
TU Kaiserslautern,
Fachbereich Mathematik,
Erwin-Schr\"odinger-Stra\ss e,
67663 Kaiserslautern, Germany, \texttt{lossen@mathematik.uni-kl.de} \and
Tel Aviv University,
School of Mathematical Sciences,
Ramat Aviv,
Tel Aviv 69978, Israel,
\texttt{shustin@math.tau.ac.il}.}

\maketitle

\begin{abstract}
In this survey, we report on progress concerning families of projective curves
with fixed number and fixed (topological or analytic) types of
singularities. We are, in particular, interested in numerical, universal and
asymptotically proper sufficient conditions to guarantee the nonemptyness,
T-smoothness and irreducibility of the variety of all projective 
curves with prescribed singularities in a fixed linear system. We also discuss
the analogous problem for hypersurfaces of arbitrary dimension with isolated
singularities, and we close with a section on open problems
and conjectures.   
\end{abstract}

\section{Introduction}

\noindent
Let $\ks$ be a topological
or analytic classification of isolated plane curve
singularities. Assume that $\Sig$ is a non-singular projective
algebraic surface over an algebraically closed
field of characteristic zero, and \mbox{$D\subset\Sig$} is an ample
divisor such that \mbox{$\dim|D|>0$} and the general member of $|D|$ is
irreducible and non-singular. The part of the discriminant,
consisting of irreducible curves on $\Sig$ with only
isolated singularities, splits into the union of {\em equisingular
families\/} ({\em ESF\/})\index{equisingular family} $\VDi$,
that is, into the varieties of irreducible curves $C\in|D|$ having exactly $r$
isolated singular points of types $S_1,\dots,S_r\in\ks$,
respectively.
If \mbox{$\Sig=\P^2$}, \mbox{$D=dH$}, $H$ a hyperplane divisor,
we simply write \mbox{$\Vdi(S_1,\dots,S_r)$} for the variety of
irreducible plane curves of degree $d$ with $r$ singularities of the prescribed
types. We focus on the following geometric problems which have been of
interest to algebraic geometers since the early 20th century but which
are still widely open in general:

\medskip\noindent
{\bf A. Existence Problem}:\index{existence of curves} Is $\VDi$
non-empty, that is, does there exist a curve \mbox{$F\in|D|$}
with the given collection of singularities\,? In particular, the
question about the minimal degree of a plane curve
having a given singularity is of special interest.

\medskip\noindent
{\bf B. T-Smoothness Problem}: If $\VDi$ is non-empty, is it
smooth and of the expected dimension (expressible via local
invariants of the singularities)\,? More precisely, let
\mbox{$C\in \VDi$} have singular points $z_1,\dots,z_r$ of types
$S_1,\dots,S_r$, respectively. We say that the family $\VDi$ is
{\em T-smooth\/}\index{tsmooth@T-smooth} at $C$ if, for every
$i=1,\dots,r$, the germ at $C$ of the family of curves
\mbox{$C'\!\in|D|$} with a singular point of type $S_i$ in a
neighbourhood of $z_i$ is smooth and has the expected dimension (to
be explained later),
and, furthermore, all these $r$ germs intersect transversally at $C$
(whence the name T-smooth).

\medskip\noindent
{\bf C. Irreducibility Problem}: Is
  $\VDi$ irreducible?

\medskip\noindent
{\bf D. Deformation Problem}: What are
the adjacency relations of
ESF in the discriminant? In other words, which simultaneous
deformations of the singularities of \mbox{$F\in \VDi$}
can be realized by a variation of $F$ in $|D|$\,?
In fact, this question is closely related to Problem B. For instance, for
\mbox{$\Sigma=\P^2$}, the
T-smoothness of \mbox{$\Vdi(S_1,\dots,S_r)$} for analytic singularities
$S_1,\dots,S_r$ is equivalent to the linear system 
$|dH|$ inducing a joint versal deformation of all singular points of
any member \mbox{$C\in \Vdi(S_1,\dots,S_r)$}. Similarly, the
T-smoothness of ESF for semiquasihomogeneous topological singularities
$S_1,\dots,S_r$ implies that the independence of simultaneous
``lower'' (w.r.t.\ the Newton diagram) deformations of the
singularities of \mbox{$C\in \Vdi(S_1,\dots,S_r)$} (see Section
\ref{sec:indep}).

\medskip\noindent
Of course, the same questions can be posed for \mbox{$V_{|D|}
  (S_1,\dots,S_r)$}, the variety of reduced (but not necessarily
irreducible) curves \mbox{$C\in|D|$} with given singularities
\mbox{$S_1,\dots,S_r\in\ks$}.

No complete solution to the above problems is known, except for the
case of plane nodal curves. However, in this overview article we
demonstrate that the approach using deformation theory and cohomology
vanishing (developed by the authors of this article over the last
ten years) enables us to obtain reasonably proper sufficient
conditions for the affirmative answers to the problems stated.
More precisely, we intend to give sufficient conditions for a positive
answer to the above questions which are
\begin{itemize}
\item {\it numerical\/}, that is, presented in the form of inequalities
relating numerical invariants of the surface, the linear system,
and the singularities,
\item {\it universal\/}, that is, applicable to curves
in any ample linear system with any
number of arbitrary singularities, and
\item {\it asymptotically proper\/} (see the definition below).
\end{itemize}
We should like to comment on the latter in more detail. Let
\mbox{$D=dD_0$} with a given divisor \mbox{$D_0\subset\Sig$}
(e.g., \mbox{$D_0=H$} if \mbox{$\Sig=\P^2$}). Then the
known general restrictions to the singular points of a curve \mbox{$C\in
\VDi$} read as upper bounds to some total singularity invariants by
a quadratic function in $d$. As an example, consider the bound
obtained by the genus formula for irreducible curves,
$$\sum_{i=1}^r\del(S_i)\,\le\,\frac{1}{2}(d^2D_0^2+dK_\Sig D_0)+1\:
,$$
where $\delta(S_i)$ is the difference between the arithmetic genus and
the geometric genus of $C$ imposed by the singularity $S_i$.
Sufficient conditions for the ``regular" properties of ESF
appear in a similar form as upper bounds to total singularity
invariants by a linear or quadratic function of the parameter $d$.
In this sense we speak of {\it linear}, or {\it quadratic}
sufficient conditions, the latter revealing the relevant
asymptotics.

Furthermore, among the quadratic sufficient conditions we
emphasize on {\it asymptotically proper\/} ones.
We say that the inequality
$$\sum_{i=1}^r a(S_i)\le f_{\Sig,D_0}(d)\ ,$$
with $a(S)$ a local invariant of singularities \mbox{$S\in\ks$},
is an {\it asymptotically proper\/}\index{asymptotically proper condition}
sufficient condition for a regular property
(such as non-emptiness, smoothness, irreducibility) of ESF, if
\begin{enumerate}
\item[(1)] for arbitrary $r,d\ge 1$
and $S_1,\dots,S_r\in\ks$ it provides the required
property of the ESF $V_{|dD_0|}^{\text{\it irr}}(S_1,\dots,S_r)$,
\item[(2)] there exists an absolute constant \mbox{$A\geq 1$}
such that, for each singularity type $S$, there is an infinite sequence
of growing integers $d,r$ and a (maybe, empty) finite collection of
singularities $\ks_d$ satisfying
$$r\cdot a(S)=A\cdot f_{\Sig,D_0}(d)+o(f_{\Sig,D_0}(d)),\quad
\sum_{S'\in\ks_d}a(S')=o(f_{\Sig,D_0}(d))\ ,$$
such that the ESF $V^{irr}_{|dD_0|}(r\cdot S,\ks_d)$ does not
have that regular property.
\end{enumerate}
If \mbox{$A=1$} in (2), we speak of an {\it asymptotically
  optimal\/}\index{asymptotically optimal condition}
sufficient condition.

In less technical terms, we say that a condition is asymptotically optimal
(resp.\ proper) if the necessary and the sufficient conditions for a regularity
property coincide (resp.\ coincide up to multiplication of the
right-hand side with a constant) if $r$ and $d$ go to infinity. 


\subsection*{Methods and Results: from Severi to Harris}

\noindent
{\it T-Smoothness.}
Already the Italian geometers \cite{Sev,Segre,Segre1} noticed
that it is possible to express the T-smoothness of the variety of
plane curves with fixed number of {\em nodes\/}\index{node} 
($A_1$-singularities)
and {\em cusps\/}\index{cusp} ($A_2$-singularities) infinitesimally.
The problem was historically called the {\em
``completeness of the characteristic linear series of complete continuous
systems''\/}  (of plane  curves with nodes and cusps).
Best known is certainly Severi's \cite{Sev} result saying that
each non-empty variety of plane nodal curves is T-smooth.
  
But for more complicated singularities (beginning with cusps) there are
examples of irreducible curves where the T-smoothness fails (see
below). That is, ether the ESF is non-smooth or its dimension exceeds
the one expected by subtracting the number of (closed) conditions imposed by
the individual singularities from the dimension $\frac{d(d+3)}{2}$ of
the variety of all curves of degree $d$. In this context, each node imposes
exactly one condition. This may be illustrated as follows: a plane curve
\mbox{$\{f=0\}$} has a node at the origin iff 
the 1-jet of $f$ vanishes ($=$\,three closed conditions) and the 2-jet
is reduced ($=$\,one open condition). Allowing the node to move (in $\C^2$)
reduces the number of closed conditions by two. Similarly, it can be
seen that each cusp imposes two conditions. 

Various sufficient conditions for
T-smoothness were found. The classical result is that the variety
\Vdnk\ of irreducible plane curves of degree $d$ with $n$ nodes and
$k$ cusps as only singularities is T-smooth, that is, smooth of
dimension \mbox{$\frac{d(d+3)}{2}-n-2k$}, if
\begin{equation}
\label{k<3d}
k\:<\:3d\,.
\end{equation}
For arbitrary singularities, several generalizations and extensions of
(\ref{k<3d}) were found \cite{GrK,GrL,Shu87,Shu91,Vas}. All of
them are of the form that the sum of certain invariants of the
singularities is bounded from above by a linear function in $d$. On
the other hand, the known restrictions for existence and T-smoothness
(and the known series of non T-smooth ESF) suggested that an
asymptotically proper sufficient condition should be quadratic in $d$ (see below).

\medskip\noindent
{\it Restrictions for the Existence.}
Various restrictions for the existence of plane curves with prescribed
singularities \mbox{$S_1,\dots,S_r$} have been found.
First, one should mention the general classical bounds
\begin{equation}
\label{genus formula}
\sum_{i=1}^r \delta(S_i) \: \leq \:
\frac{(d\!\!\;-\!\!\;1)(d\!\!\;-\!\!\;2)}{2}\,,
\end{equation}
(for irreducible curves), resulting from the genus formula, respectively
$$ \sum_{i=1}^r \mu(S_i)  \: \leq \: (d\!\!\;-\!\!\;1)^2, $$
resulting from the intersection of two generic
polars and B{\'e}zout's theorem.
If $C$ has only nodes and cusps as singularities then the
Pl\"ucker formulas give (among
others) the necessary conditions
\begin{eqnarray*}
2 \cdot \# (\text{nodes\/}) + 3 \cdot \# (\text{cusps\/}) & \leq &
d^2\!\!\:-d-2\,, \\
6 \cdot \# (\text{nodes\/}) + 8 \cdot \# (\text{cusps\/}) & \leq & 3d^2
\!\!\:-6d\,.
\end{eqnarray*}
By applying the log-Miyaoka inequality, F.\ Sakai \cite{Sak}
obtained the necessary condition
$$ \sum_{i=1}^r \mu(S_i) \: < \: \frac{2\nu}{2\nu+1} \cdot \left(d^2\! -
\frac{3}{2}\!\: d \right) \,,$$
where $\nu$ denotes the maximum of the multiplicities $\mt S_i$,
\mbox{$i=1,\dots,r$}.
If \mbox{$S_1, \dots,S_r$} are ADE-singularities then
$$ \sum_{i=1}^r \mu(S_i) \: < \: \left\{
\renewcommand{\arraystretch}{2.0}
\begin{array}{cl}
\displaystyle{
\frac{3}{4}\!\:d^2\! -\frac{3}{2} \!\: d +2} & \text{ if $d$ is
even}\,,\\
\displaystyle{\frac{3}{4}\!\:d^2\! - d +\frac{1}{4}} & \text{ if $d$
  is odd}\,,\\
\end{array}
\right.
$$
is necessary for the existence of a plane curve with $r$ singularities of types
\mbox{$S_1, \dots,S_r$} (cf.\ \cite{HF,Ivi}, resp.\ \cite{Sak}). Further
necessary conditions can be obtained, for instance, by applying the
semicontinuity of the singularity spectrum (see \cite{Var1}).

\medskip\noindent
{\it Methods of Construction.}
The first method is to construct (somehow) a curve of the given degree
which is degenerate with respect to the required curve, and then to
deform it in order
to obtain the prescribed singularities.

\begin{figure}[ht]

\bigskip
\begin{center}
\input{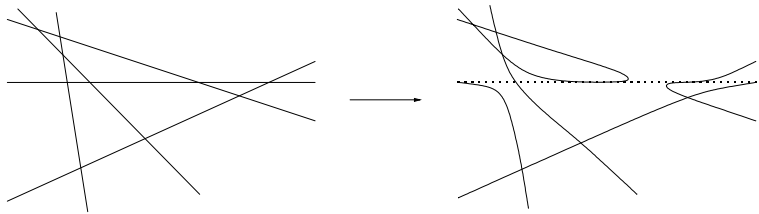} 
\caption{\label{fig:lines} Constructing irreducible nodal curves.}
\end{center}
\end{figure}

For instance, Severi \cite{Sev} showed that singular
points of a nodal curve,
irreducible or not, can be smoothed, or preserved, independently.
Hence, starting with the union of generic straight lines in the
projective plane and smoothing suitable
intersection points, one obtains irreducible curves with any
prescribed number of nodes, allowed by the genus bound (\ref{genus formula}), see
Fig.\ \ref{fig:lines}.

\noindent
Attempts to extend this construction to other singularities give curves
with a number of singularities
bounded from above by a linear function in the degree $d$
(see, for example, \cite{GrM} for curves with nodes,
cusps and ordinary triple points), because of the very restrictive
requirement of the independence of deformations.

The second method consists of a construction especially adapted to
the given degree and given collection of singularities. It may be
based on a sequence of birational transformations of the plane
applied to a more or less simple initial curve in order to obtain
the required curve.
Or it may consist of an invention of a polynomial defining the
required curve.
This is illustrated by constructions of singular curves
of small degrees as, for instance, in \cite{Wl1}, \cite{Wl2},
or by Hirano's \cite{H} construction of cuspidal curves, which
led to a series of irreducible
cuspidal curves of degree \mbox{$d=2\!\!\:\cdot \!\!\: 3^k$},
\mbox{$k\geq 1$} with precisely \mbox{$9(9^k\!\!\:-\!\!\:1)/8$}
cusps. Note that in this case the number of conditions
imposed by the cusps is \mbox{$d^2/16+O(d)$} more than the dimension of the space of
curves of degree $d$. 

Two main difficulties do not allow to apply this
approach to a wide class of degrees and singularities:
\begin{itemize}
\item for any
new degree or singularity one has to invent a new construction,
\item even if one has constructed a curve with many singularities, it
  is hard to check that these singular points can be
smoothed independently or, at least, that any intermediate number of singularities
can be realized (for instance, for Hirano's examples the latter can
hardly be expected).
\end{itemize}

\medskip\noindent
{\it Irreducibility.}
Severi \cite{Sev} claimed that all non-empty ESF
\mbox{$\Vdi(n\!\!\:\cdot\!\!\: A_1)$}
of irreducible plane nodal curves are not only T-smooth but also
irreducible. However, as was realized later, his proof of the
irreducibility was incomplete, and the problem has become known later
as {\em ``Severi's conjecture''}.

For many years, algebraic geometers tried to solve this problem
without much progress. A major step forward was made by Fulton
\cite{Fulton} and Deligne \cite{Deligne}, showing that the fundamental
group of the
complement of each irreducible nodal curve is Abelian (which is a
necessary condition for \mbox{$\Vdi(n\cdot A_1)$} to be irreducible). Finally, the
problem was settled by Harris \cite{Harris}. He gave a rigorous proof
for the irreducibility of the varieties
\mbox{$\Vdi(n\!\!\:\cdot\!\!\: A_1)$}, by inventing a new
specialization method and by using the irreducibility of the moduli
space of curves of a given genus.
Since nodal curves form an open dense subset
of each {\em Severi variety\/}\index{Severi variety} \cite{Alb,Nob1},
that is, of the variety of all irreducible plane curves of a
given degree and genus, Harris' theorem extends to
all Severi varieties as well. Later Ran \cite{Ran2} and
Treger \cite{Tr}
gave different proofs. Then Ran \cite{Ran} generalized Harris'
theorem to ESF \mbox{$\Vdi(O_m,n\!\!\:\cdot\!\!\: A_1)$} with one
ordinary singularity $O_m$ (of some order \mbox{$m\geq 1$}) and any
number of nodes \mbox{$n\le(d-1)(d-2)/2-\del(O_m)$}. Kang
\cite{Kang,Kang2} obtained 
the irreducibility of the families \mbox{$\Vdnk$} for
$$k\le 3\,,\quad\text{or}\quad \frac{d^2\! - 4d +1}{2}\le n \le
\frac{d^2\! -3d + 2}{2}\ .$$
The main idea of the proofs consists of using moduli spaces
of curves, special degenerations of nodal curves, or degenerations of
rational surfaces. In any case the proofs heavily rely on the
independence of simultaneous deformations
of nodes for plane curves, which (in general) does not hold for more
complicated singularities, even for cusps.

\medskip
\noindent
{\it Examples of Obstructed and Reducible ESF.}
Already, Segre \cite{Segre1,Tan84}
constructed a series of irreducible plane curves such
that the corresponding germs of ESF are non-T-smooth:
\mbox{$V_{6m}^{\text{\it irr}} (6m^2 \!\cdot\!\!\:A_2)$}, \mbox{$m\geq
  3$}. Similar examples are given in \cite{Shu94}.
However, in these examples $\Vd$ is smooth (but of bigger dimension than the
expected one).
In 1987, Luengo \cite{Luengo1} provided  the first examples of curves $C$ such
that the corresponding ESF is {\em of expected dimension, but
  non-smooth}, e.g.\ \mbox{$V_{9}^{\text{\it irr}}(A_{35})$}.

Concerning reducible ESF, there is mainly one classical example due to
Zariski \cite{Zar}: \mbox{$V^{\text{\it irr}}_{6} (\!\:6 \!\!\:\cdot \!\!\: A_2)$}.
More precisely, Zariski shows that there exist exactly two
components, both being T-smooth. One component consists of cuspidal
sextics $C$ whose singularities lie on a conic.
For such curves, Zariski computed the fundamental
group \mbox{$\pi_1(\P^2\setminus C)$} to be the {\em non-Abelian\/} group
\mbox{$\Z_2 \ast \Z_3$}. On the other hand, he showed that there exist
curves $C$ whose singularities do not lie all on a conic and
whose complement in $\P^2$
has an {\em Abelian\/} fundamental group. In particular, these curves
cannot be obtained by a deformation from curves in the
first component.

Actually, this example belongs to a series \cite{Shu94}:
\mbox{$V^{\text{\it irr}}_{6p} (\!\:6p^2
\!\!\:\cdot \!\!\: A_2)$}, \mbox{$p\geq 2$}, is reducible.
More precisely, for \mbox{$p\geq 3$} there exist
components with different dimensions. For \mbox{$p=2$} there exist two
different T-smooth components, as in Zariski's example.

\section{Geometry of ESF in Terms of Cohomology}
\label{sgdt1}
\setcounter{equation}{0}

\noindent
In this section, we relate the ESFs $\VD$ to strata of a Hilbert
scheme representing a deformation functor. This allows us to deduce
geometric properties such as the T-smoothness or the irreducibility from the
vanishing of the first cohomology group for the ideal sheaves of
appropriately chosen zero-dimensional subschemes of the surface $\Sigma$.

Let $T$ be a complex space, then by a {\em family of reduced
  (irreducible) curves on $\Sigma$ over}\index{family of curves} $T$ we 
mean a commutative diagram
$$
\xymatrix@C=6pt@R=10pt{ {\sC} \:\:\ar@{^{(}->}[rr]^-{j}
 \ar[rd]_\varphi
 & & \:\Sigma \times T \ar[dl]^{\pr}\\
&  T }
$$
where $\varphi$ is a proper and flat morphism such that all
 fibres $\sC$$_t:=\varphi^{-1}(t)$, \mbox{$t \in T$}, are
reduced (irreducible) curves on $\Sigma$, \mbox{$j : \sC
\hookrightarrow$$\, \Sigma \times T$} is a closed embedding and $\pr$
denotes the natural projection.

A {\em family with sections\/}
is a diagram as above, together with sections
\mbox{$\sigma_1,\dots,\sigma_r:T\to \sC$} of $\varphi$.

To a family of reduced plane curves and a point \mbox{$t_0\in T$}
we can associate, in a functorial way, the deformation
\mbox{$ ({\sC},z_1){\scriptstyle\amalg} \ldots
  {\scriptstyle\amalg}({\sC},z_r) \rightarrow (T,t_0)$} 
of the multigerm \mbox{$(C,  \Sing\, C)=
\coprod_i (C,z_i)$} over the germ $(T,t_0)$, where
\mbox{$C={\sC}_{t_0}$} is the fibre over $t_0$. 
Having a family with sections
\mbox{$\sigma_1,\dots,\sigma_r$}, \mbox{$\sigma_i(t_0)=z_i$}, we
obtain in the same way a deformation  of \mbox{$\coprod_i (C,z_i)$}
over $(T,t_0)$ with  sections.

A family \mbox{$\sC \hookrightarrow$$\, \Sigma\times T \to T$} of
reduced  curves (with sections) is called {\em
  equianalytic}\index{equianalytic family}
(along the sections) if, for each \mbox{$t\in T$}, the induced  
deformation of the multigerm \mbox{$(\sC_t, \Sing\, \sC_t)$} is 
isomorphic (isomorphic as deformation with section) to the trivial
deformation (along the trivial sections).
It is called {\em equi\-singular}\index{equisingular family}
(along the sections) if, for each \mbox{$t\in T$}, the induced 
deformation of the multigerm \mbox{$(\sC_t, \Sing\, \sC_t)$} is 
isomorphic (isomorphic as deformation with section) to 
an equisingular deformation along the trivial sections. In other words,
a family with sections is equianalytic (resp.\ equisingular) if the
analytic (resp.\ topological) type of the fibre $\sC_t$ does not change along the
sections (cf.\ \cite{Wahl1,GLSBook1}). 
Here, the {\em analytic type\/}\index{analytic type} of a reduced plane curve
singularity $(C,z)$ is given by the isomorphism class of its analytic
local ring $\ko_{C,z}$, while the {\em (embedded) topological
  type\/}\index{topological type} is given by
the Puiseux pairs of its branches and their mutual intersection
multiplicities or, alternatively, by the system of multiplicity
sequences (see \cite{BrK}).

The {\em Hilbert functor\/}\index{Hilbert functor}
$\Hilb_{\Sigma}$ on the category of complex
spaces defined by
$$ \Hilb_{\Sigma}(T):=\{\!\; {\sC} \hookrightarrow \Sigma\times T
\to T, \mbox{ family of reduced curves over $T$} \!\;\}$$
is known to be representable by a complex space $\tHilb_{\Sigma}$
(see \cite{Grothendieck} for algebraic varieties and \cite{Douady} for
arbitrary complex spaces).
Moreover, the universal family of reduced curves on $\Sigma$ ``breaks up'' into
strata with constant Hilbert polynomials, more precisely,
\mbox{$  \tHilb_{\Sigma}\! = \coprod_{h\in \K\!\:[z]}
\tHilb_{\Sigma}^h $}
where the $\tHilb_{\Sigma}^h$
are (unions of) connected components of $\tHilb_{\Sigma}$ whose points
correspond to curves on $\Sigma$ with the fixed Hilbert polynomial $h$.

Let \mbox{$V_h(S_1,\dots,S_r)\subset \tHilb_{\Sigma}^h$} denote the
locally closed subspace ({\it equisingular stratum\/}) of reduced
curves with Hilbert polynomial $h$ having precisely $r$ singularities of types
\mbox{$S_1,\dots,S_r$} (\cite{GrL}). Further, let
{$\Vhi(S_1,\dots,S_r)\subset
  V_h(S_1,\dots,S_r)$} denote the open subspace parametrizing
    irredu\-cible curves.

This notion of (not necessarily reduced) equisingular strata in the Hilbert scheme
$\tHilb_{\Sigma}$ is closely
related to the ESFs $\VD$ considered before: if \mbox{$U \subset
  |C|=\P\bigl(H^0(\ko_{\Sigma}(C))\bigr)$}
denotes the open subspace corresponding to reduced curves, then there
exists a unique morphism \mbox{$U\hookrightarrow
\tHilb_{\Sigma}^h$} which on the tangent level corresponds to
\mbox{$ H^0\bigl(\ko_{\Sigma}(C)\bigr) \big/
H^0\bigl(\ko_{\Sigma}\bigr) \hookrightarrow
H^0\bigl(\ko_{C}(C)\bigr)$}. Via this morphism, we may consider $U$ as
a locally closed subscheme of 
\mbox{$\tHilb_{\Sigma}^h$}. In particular, for a regular surface
(that is, a surface satisfying \mbox{$H^1(\ko_\Sigma)=0$}) the above
injection is an isomorphism and $U$ is an
open subscheme of \mbox{$\tHilb_{\Sigma}^h$} (see \cite{GrL1,GLSBook1} for
details).

In the following, we give a geometric interpretation of the zeroth
  and first cohomology of the ideal sheaves of certain
  zero-dimensional schemes. We write \mbox{$C\in
  V_h(S_1,\dots,S_r)$} to denote either the point in
$V_h(S_1,\dots,S_r)$ or the curve corresponding to the point, that is,
the corresponding fibre of the universal family
${\sU}_h$ over the Hilbert scheme. Consider the map
\begin{equation}
\label{Phid}
\Phi_h:  V_h(S_1,\dots,S_r) \lra \mbox{Sym}^r\Sigma\,, \quad C
\longmapsto (z_1\!+\!\ldots\!+\!z_r)\,,
\end{equation}
where $\mbox{Sym}^r\Sigma$
is the $r$-fold symmetric product of $\Sigma$ and where
\mbox{$(z_1\!+\!\ldots \!+\!z_r)$} is the non-ordered tuple of the
singularities of $C$.
Since each equisingular, in particular each equianalytic, deformation of
a germ admits a unique singular section (cf.\ \cite{Tei}), the
universal family over $V_h(S_1,\dots,S_r)$,
$${\sU}_h(S_1,\dots,S_r)\hookrightarrow \Sigma \times V_h(S_1,\dots,S_r)
\to V_h(S_1,\dots,S_r),$$ admits, locally at $C$, $r$ singular
sections. Composing these sections with the projections to $\Sigma$
gives a local description of the map $\Phi_h$ and shows in particular
that $\Phi_h$ is a well-defined morphism, even if
\mbox{$V_h(S_1,\dots,S_r)$} is not reduced.

We denote by \mbox{$\Vhf(S_1,\dots,S_r)$} the complex space consisting of
the disjoint union of the fibres of $\Phi_h$. Thus, each
connected component of \mbox{$\Vhf(S_1,\dots,S_r)$} consists of curves
with fixed positions of the singularities in $\Sigma$.

It follows from the universal property of
$V_h(S_1,\dots,S_r)$ and from the above construction that
\mbox{$\Vhf (S_1,\dots,S_r)$}, together with the
induced universal family on each fibre, represents the
functor of equianalytic, resp.\ equisingular, families of given
types $S_1,\dots,S_r$ along trivial sections.

Before formulating the main proposition relating the vanishing of cohomology to
geometric properties of ESF, we introduce some notation: we write
$V_h$ to denote \mbox{$V_{h}(S_1,\dots,S_r)$}, resp.\
\mbox{$\Vhf(S_1,\dots,S_r)$}, and $V_{|C|}$ to denote
\mbox{$\VC$}, resp.\ \mbox{$\VCf(S_1,\dots,S_r)$}.
Moreover, we write 
$$Z'(C)=\bigcup_{z\in \Sing(C)} Z'(C,z)$$ to denote
one of the $0$-dimensional schemes $\Zea (C), \Zeaf (C), \Zes (C),
\Zesf (C)$, where
\begin{itemize}
\itemsep3pt
\item $\Zea (C,z)$ is defined by the {\em Tjurina ideal\/}, that is,
  the ideal generated by a local equation \mbox{$f\in \C\{u,v\}$} for
  \mbox{$(C,z)\subset (\Sigma,z)$} and its partial derivatives;
\item $\Zeaf (C,z)$ is defined by the ideal
\mbox{$\Ieaf(f)=\langle f \rangle +
\langle u,v\rangle \!\cdot \!\langle \tfrac{\partial f}{\partial u},\tfrac{\partial f}{\partial v}\rangle$};
\item $\Zes (C,z)$ is defined by the {\em equisingularity ideal\/}
  \mbox{$\Ies(f)\subset \C\{u,v\}$} as introduced by Wahl \cite{Wahl1};
\item $\Zesf (C,z)$ is defined by the ideal
$$
\Iesf (f):=\left\{ g\!\!\;\in\!\!\; \Ies(f) \,\bigg|\:
\renewcommand{\arraystretch}{0.9}
  \begin{array}{c}
f +\varepsilon g \text{ defines an
  equisingular deformation}\\
\text{of $(C,0)$ along the trivial section over $T_\varepsilon$}
\end{array}
\!\!\:\right\}
$$
(here, $T_\varepsilon$ denotes the fat point $(\{0\},\C[\varepsilon]/\langle
\varepsilon^2])$).
\end{itemize}

\noindent
We write \mbox{$\kj_{Z'(C)/C}$}, resp.\  $\kj_{Z'(C)/\Sigma}$,
    to denote the ideal sheaf of $Z'(C)$ in $\ko_C$, resp.\ in
    $\ko_{\Sigma}$, and \mbox{$\kj(C):=\kj\otimes_{\ko_{\Sigma}}
      \ko_{\Sigma}(C)$}. Moreover, we write $\deg Z'(C)$ for the
    degree of $Z'(C)$ as a projective variety, that is, 
$$\deg
    Z'(C)=\dim_{\C} \ko_{\Sigma}/\kj_{Z'(C)/\Sigma}\,.$$

\begin{proposition}[{\cite[Prop.\:2.6]{GrL1}}]
\label{H0H1prop}
Let \mbox{$C\subset \Sigma$} be a reduced curve with Hilbert
polynomial $h$ and precisely $r$ singularities
\mbox{$z_1,\dots,z_r$} of analytic or topological types
\mbox{$S_1,\dots,S_r$}.
\begin{description}
\itemsep3pt
\item[(a)] The Zariski tangent space of $V_h$ at $C$ is
  \mbox{$H^0\bigl(\kj_{Z'(C)/C}(C)\bigr)$}, while
 the Za\-ris\-ki tangent space of $V_{|C|}$ at $C$ is
  \mbox{$H^0\bigl(\kj_{Z'(C)/\Sigma}(C)\bigr)
\big/H^0\bigl(\ko_{\Sigma}\bigr)$}.
\item[(b)] \mbox{$h^0\bigl(\kj_{Z'(C)/C}(C)\bigr)-
h^1\bigl(\kj_{Z'(C)/C}(C)\bigr) \:\leq \: \dim (V_h,C)\leq
h^0\bigl(\kj_{Z'(C)/C}(C)\bigr)$}.
\item[(c1)] If \mbox{$H^1\bigl(\kj_{Z'(C)/C}(C)\bigr)=0$} then $V_h$ is
T-smooth at $C$, that is,
smooth of the expected dimension $h^0\bigl(\ko_C(C)\bigr)-\deg
Z'(C) = C^2+1-p_a(C) - \deg Z'(C)\,.$
\item[(c2)]  If \mbox{$H^1\bigl(\kj_{Z'(C)/\Sigma}(C)\bigr)=0$} then
$V_{|C|}$ is T-smooth at $C$, that is, smooth of the expected dimension
  \mbox{$h^0\bigl(\ko_{\Sigma}(C)\bigr)-1-\deg Z'(C)$}.
\item[(d)]
If \mbox{$H^1\bigl(\kj_{\Zea (C)/C}(C)\bigr)=0$} then
the natural morphism of germs
$$ \bigl( \tHilb^h_\Sigma,\,C\bigr) \lra
\textstyle{\prod\limits_{i=1}^r}
\text{Def}\,(C,z_i)$$
is smooth of fibre dimension
\mbox{$h^0\bigl(\kj_{\Zea(C)/C}(C)\bigr)$}.
Here, \mbox{$\prod_{i=1}^r \text{Def}\,(C,z_i)$} is the Cartesian
product of the base spaces of the semiuniversal deformations of the
germs $(C,z_i)$.
\item[(e)] Write \mbox{$\Zf(C)$} for
$\Zeaf(C)$, respectively  $\Zesf(C)$. Then  the vanishing of
\mbox{$H^1\bigl(\kj_{\Zf(C)/C}(C)\bigr)$}
implies that the morphism of germs
$$\Phi_h: \bigl(V_h(S_1,\dots,S_r),C\bigr) \rightarrow
 \bigl(\mathrm{Sym}^r\Sigma,(z_1\!+\!\ldots\!+\!z_r)\bigr) $$
is smooth of fibre dimension
  \mbox{$h^0\bigl(\kj_{\Zf(C)/C}(C)\bigr)$}.
\end{description}
\end{proposition}

\noindent
We reformulate and strengthen Proposition \ref{H0H1prop} in the
case of plane curves, which is of special interest. Of  course, since
\mbox{$h^1\bigl(\ko_{\P^2}\bigr)=h^2\bigl(\ko_{\P^2}\bigr)=0$},  there
is no difference whether we consider the curves in the linear
system $|dH|$, $H$ the hyperplane divisor, or curves with fixed
Hilbert polynomial \mbox{$h(z)=dz-(d^2\!-\!\!\:3d)/2$}. We denote the
corresponding varieties by \mbox{$V_d= V_d(S_1,\dots,S_r)$},
respectively by \mbox{$V_{d,\mbox{\scriptsize{fix}}}(S_1,\dots,S_r)$}. Using
the above notation, we obtain:

\begin{proposition}[{\cite[Prop.\:2.8]{GrL1}}]
\label{H0H1}
Let \mbox{$C\subset \P^2$} be a reduced curve of degree
$d$ with precisely $r$ singularities \mbox{$z_1,\dots,z_r$} of
analytic or topological types \mbox{$S_1,\dots,S_r$}.
\begin{enumerate}
\itemsep3pt
\item[(a)] \mbox{$H^0\bigl(\kj_{Z'(C)/\P^2}(d)\bigr)\big/ H^0(\ko_{\P^2})$}
is isomorphic to the Zariski tangent space of $V_d$ at $C$.
\item[(b)] \mbox{$h^0\bigl(\kj_{Z'(C)/\P^2}(d)\bigr)-
h^1\bigl(\kj_{Z'(C)/\P^2}(d)\bigr)-1 \:\leq \: \dim (V_d,C)$}\\
\phantom{$h^0\bigl(\kj_{Z'(C)/\P^2}(d)\bigr)-
h^1\bigl(\kj_{Z'(C)/\P^2}(d)\bigr)-1$} \mbox{$\:\leq \:
h^0\bigl(\kj_{Z'(C)/\P^2}(d)\bigr)-1$}.
\item[(c)]\mbox{$H^1\bigl(\kj_{Z'(C)/\P^2}(d)\bigr)=0$} iff $V_d$ is
T-smooth at $C$, that is,
smooth of the expected dimension \mbox{$d(d\!\!\:+\!\!\:3)/2 - \deg
  Z'(C)$}.
\item[(d)]\mbox{$H^1\bigl(\kj_{\Zea (C)/\P^2}(d)\bigr)=0$} iff the natural
  morphism of germs
$$\left( \P\bigl(H^0\bigl(\ko_{\P^2}(d)\bigr)\bigr),\,C\right) \to
\textstyle{\prod_{i=1}^r} \text{Def}\,(C,z_i)$$
is smooth (hence surjective)
of fibre dimension \mbox{$h^0\bigl(\kj_{\Zea(C)/\P^2}(d)\bigr)-1$}.
\item[(e)] Write \mbox{$\Zf(C)$} for
$\Zeaf(C)$, respectively  $\Zesf(C)$.
Then
$H^1\bigl(\kj_{\Zf(C)/\P^2}(d)\bigr)$
vanishes iff the morphism of germs
$$\Phi_d: \bigl(V_d(S_1,\dots,S_r),C\bigr) \rightarrow
 \bigl(\mathrm{Sym}^r\P^2,(z_1\!+\!\ldots\!+\!z_r)\bigr) $$
is smooth of fibre dimension
  \mbox{$h^0\bigl(\kj_{\Zf(C)/\P^2}(d)\bigr)-1$}.
In particular, the vanishing of
\mbox{$H^1\bigl(\kj_{\Zf(C)/\P^2}(d)\bigr)$} implies that
arbitrarily close to $C$ there are curves in $V_d(S_1,\dots,S_r)$
whose singularities are in general position in $\P^2$.
\end{enumerate}
\end{proposition}

\section{T-Smoothness} \label{sec:smoothness}

To show the T-smoothness of $\VC$ it suffices to show, according to Proposition
\ref{H0H1prop}\:(c2), that  
\mbox{$H^1\bigl(\kj_{\Zea(C)/\Sigma}(C)\bigr)=0$} (in the case of analytic
types), respectively  \mbox{$H^1\bigl(\kj_{\Zes(C)/\Sigma}(C)\bigr)=0$} (in
the case of topological types). Note that for \mbox{$\Sigma=\P^2$},
these conditions are even equivalent to the T-smoothness of $V_d(S_1,\dots,S_r)$
by Proposition \ref{H0H1}\:(c).

\subsection{ESF of Plane Curves}

The classical approach to the $H^1$-vanishing problem (based on
Riemann-Roch and Serre duality) leads to sufficient conditions for the
T-smoothness of ESF of plane curves such as the $3d$-condition
(\ref{k<3d}) and its extensions mentioned above. 

In the papers \cite{GLS97,GLS98}, we applied two different approaches to the
$H^1$-vanishing problem, based on the Reider-Bogomolov theory
of unstable rank 2 vector bundles (see also \cite{ChS}), respectively
on the Castelnuovo function of the ideal sheaf of a zero-dimensional
scheme (see also \cite{Davis,Barkats}). Both approaches lead to
quadratic sufficient conditions for the T-smoothness of ESF of plane
curves. Combining both approaches, we obtain:

\begin{theorem}[\cite{GLS98,GLS01}]
\label{Smoothness Theorem}
Let \mbox{$C\subset \P^2$} be an irreducible curve of degree
\mbox{$d>5$} having $r$ singularities \mbox{$z_1,\dots,z_r$} of
topological (respectively analytic) types $S_1,\dots,S_r$ as
its only singularities. Then \mbox{$\Vdi(S_1,\dots, S_r)$}
is T-smooth at $C$ if
\begin{equation}
\label{Smoothness condition}
\sum\limits_{i=1}^r \gamma'(C,z_i) \leq (d+3)^2,
\end{equation}
\mbox{$\gamma'(C,z_i)=\ges(C,z_i)$} for $S_i$ a topological
type, resp.\ \mbox{$\gamma'(C,z_i)=\gea(C,z_i)$} for $S_i$ an
analytic type.
\end{theorem}

\noindent
Here, $\ges$ and $\gea$ are new analytic invariants of singularities
which are defined as follows (see \cite{KeilenLossen} for a thorough discussion):

\medskip
\noindent
{\bf The $\bgamma$-invariant.} 
Let \mbox{$f\in \C\{u,v\}$} be a reduced power series, and let {$I\subset
  \langle u,v\rangle \subset \C\{u,v\}$} be an ideal  
containing the Tjurina ideal \mbox{$\Iea (f)=\langle f,\frac{\partial
  f}{\partial u},\frac{\partial f}{\partial v}\rangle$}. For each
\mbox{$g\in \langle u,v\rangle \subset \C\{x,y\}$}, we introduce 
$ \Delta (f,g;I)$ as the minimum among  \mbox{$\dim_{\C}
 \C\{u,v\}/\langle I, g\rangle $} and  \mbox{$i\!\:(f,g)-
 \dim_{\C} \C\{u,v\}/\langle I, g\rangle$}  (where
\mbox{$i\!\:(f,g)$} denotes the intersection multiplicity of $f$ and $g$,
\mbox{$i\!\:(f,g)=\dim_{\C} \C\{u,v\}/\langle f,g\rangle$}). 
By \cite[Lemma 4.1]{Shu97}, this minimum is at least $1$ so that we
may define
$$
 \gamma\!\; \bigl(f;I\bigr) := \max_{g\in \langle x,y\rangle}
 \left\{ \frac{(\dim_{\C} \C\{u,v\}/\langle I, g\rangle+\Delta
     (f,g;I))^2}{\Delta (f,g;I)}
\right\}\,,
$$
and, finally,
\begin{eqnarray*}
\gea(f)&:=& \max \bigl\{\gamma
    \!\;\bigl(f;\!\:I\bigr)\,\big|\, I \supset \Iea(f)\text{ a complete
    intersection ideal}\bigr\},\\
\ges(f)&:=& \max \bigl\{\gamma
    \!\;\bigl(f;\!\:I\bigr)\,\big|\, I \supset \Ies(f)\text{ a complete
    intersection ideal}\bigr\}.
\end{eqnarray*}

\medskip\noindent
Note that \mbox{$\gamma'(f)\leq \bigl(
\tci'(f)+1\bigr)^2$}, where $\tci'(f)$ stands for one of
\begin{eqnarray*}
 \tci(f) &:=& \max\bigl\{\!\;\dim_\C \C\{u,v\}/I \,\big|\: I\supset
 \Iea(f) \text{ a complete intersection ideal}\,\bigr\}\\ &\leq& \tau(f)\,, \\
 \tesci(f) &:=& \max\bigl\{\!\;\dim_\C \C\{u,v\}/I \,\big|\: I\supset
 \Ies(f) \text{ a complete intersection ideal}\,\bigr\}\\ &\leq& \tes(f)\,.
\end{eqnarray*}
Here, \mbox{$\tau(f)=\dim_\C
\C\{u,v\}/\Iea(f)$} is the {\em Tjurina number\/} of $f$ and
{$\tes(f)=\dim_\C \C\{u,v\}/\Ies(f)$} is the codimension of
the $\mu$-constant stratum in the semiuniversal deformation of $f$
(see \cite[Lemma 4.2]{GLS98}).

In general, these bounds for the $\gamma$-invariant are far from being
sharp. For instance, if $f$ defines an {\em ordinary singularity of order\/}
\mbox{$m\geq 3$} (that is, the $m$-jet of $f$ is a reduced homogeneous polynomial of
degree $m$), then \mbox{$\ges(f)=2m^2$} while \mbox{$(\tes(f)+1)^2
  =\frac{1}{2}m^4 +O(m)$}.  See \cite{KeilenLossen}
for the case of semiquasihomogeneous singularities.

For ESF of irreducible curves with nodes and cusps, respectively for
ESF of irreducible curves with ordinary singularities, we obtain:

\begin{corollary}
\label{Corollary 3.4}
\mbox{$V_d(n\!\!\:\cdot \!\!\:A_1,\,k\!\!\:\cdot\!\!\: A_2)$}
is T-smooth or empty if
\begin{equation}
4n+9k \: \leq \:(d+3)^2\,. \label{en1}\end{equation}
\end{corollary}

\begin{corollary}
If $S_1,\dots,S_r$ are ordinary singularities of order
$m_1,\dots,m_r$, then
\mbox{$V_d (S_1,\dots,S_r)$} is  T-smooth or empty if
\begin{equation}
4\cdot \#(\mbox{nodes\/}) +
\sum\limits_{m_i\geq 3} 2\cdot m_i^2 \: \leq \: (d+3)^2.
\label{en11}\end{equation}
\end{corollary}

\noindent
In particular, it follows that Theorem \ref{Smoothness Theorem}
is asymptotically proper for
ordinary singularities, since the inequality
$$\sum_{i=1}^r m_i(m_i\!-\!\!\;1) \: \leq \:
(d\!\!\;-\!\!\;1)(d\!\!\;-\!\!\;2)$$
is necessary for the existence of an irreducible curve with ordinary
singularities of multiplicities \mbox{$m_1,\dots,m_r$}.
More generally, by constructing series of ESF where the
T-smoothness fails (see \cite{Shu97,GLS97,LossenPrep}), we 
proved the asymptotic properness of condition \eqref{Smoothness
  condition} in Theorem \ref{Smoothness Theorem} for the case of
semiquasihomogeneous singularities.

\subsection{ESF of Curves on Smooth Algebraic Surfaces}
\label{sec:smo_other}

\noindent
Let $\Sigma$ be a smooth projective surface and $D_0$ an effective
divisor on $\Sigma$. In this situation, for lack of a generalization
of the Castelnuovo function approach, so far only the Reider-Bogomolov
approach leads to a quadratic sufficient condition for the
T-smoothness of (non-empty) ESF \mbox{$V_{|dD_0|}(S_1,\dots,S_r)$}.
Set
$$ A(\Sigma,D_0):=\frac{(D_0.K_\Sigma)^2\!-\!\!\:D_0^2K_\Sigma^2}{4}\,,$$
where $K_\Sigma$ is the canonical divisor on $\Sigma$.
By the Hodge index theorem, this is a non-negative number if $D_0$ or
$K_{\Sigma}$ is ample.

\begin{theorem}[{\cite{GLS97,GLSBook2}}]
\label{Theorem 1}
Let \mbox{$C\subset \Sigma$} be
an irreducible curve with precisely $r$ singular points
\mbox{$z_1,\dots,z_r$} of topological or analytic types
\mbox{$S_1,\dots, S_r$}, ordered such that 
\mbox{$\tci'(C,z_1)\geq \ldots \geq \tci'(C,z_r)$}. Assume that $C$ and
\mbox{$C\!\!\:-\!\!\:K_{\Sigma}$}
are ample and that \mbox{$C^2\geq \max \,\{K_{\Sigma}^2, A(\Sigma,C)\}$}. If
\begin{equation}
\label{smoot1.0}
\sum_{i=1}^r \tci'(C,z_i) \,< \,\frac{(C\!\!\:-\!\!\:K_\Sigma)^2}{4}
\end{equation}
and, for each \mbox{$1\leq s\leq r$},
\begin{equation}
\label{smoot1.1}
\Biggl( \sum_{i=1}^s
\bigl( \tci'(C,z_i)\!\!\:+\!\!\:1\bigr)\Biggr)^2 \;< \;
\sum_{i=1}^s \Bigl(C^2\!-
C.K_\Sigma \bigl( \tci'(C,z_i)\!\!\;+\!\!\;1\bigr)\Bigr)- A(\Sigma,C)\,,
\end{equation}
then \mbox{$V_{|C|}(S_1,\dots,S_r)$}
is T-smooth at $C$. Here, $\tci'(S_i)$ stands for $\tesci(S_i)$ if $S_i$ is a
topological type and for $\teaci(S_i)$ if $S_i$ is an analytic type.
\end{theorem}

\noindent
Note that, in general, the conditions (\ref{smoot1.1}) are not {\em
  quadratic\/} sufficient conditions in the above sense. But in many
cases they are. For instance, if $-K_\Sigma$ is nef, then by 
applying the Cauchy inequality we deduce: 

\begin{corollary}\label{cor:Theorem 1}
Let $\Sigma$ be a smooth projective surface with $-K_\Sigma$ nef,
$D_0$ an ample divisor on $\Sigma$ such that \mbox{$D_0^2\geq
  A(\Sigma,D_0)$}, and \mbox{$d>0$} such that 
\mbox{$d^2\geq K_\Sigma^2/D_0^2$}. If
\begin{equation}
\label{smoot1.1'}
\sum_{i=1}^r
\bigl( \tci'(S_i)\!\!\:+\!\!\:1\bigr)^2 \;< \;
\bigl(D_0^2-A(\Sigma,D_0)\bigr) \cdot d^2 - 
2(D_0.K_\Sigma)\cdot d\,, 
\end{equation}
then \mbox{$V_{|dD_0|}^{\text{\it irr}}(S_1,\dots,S_r)$}
is T-smooth or empty.
\end{corollary}

\noindent
For ESF of curves with nodes and cusps, respectively for ESF of
curves with ordinary singularities, the obvious estimates for $\tau'$
allow us to deduce the following corollaries:

\begin{corollary}\label{cn1} With the assumptions of Corollary \ref{cor:Theorem 1}, let
\begin{equation}
\label{smoot1.1'nc}
4n+9k \;< \; \bigl(D_0^2- A(\Sigma,D_0)\bigr)\cdot d^2 -
2(D_0.K_\Sigma) \cdot d 
\end{equation}
then \mbox{$V_{|dD_0|}^\text{\it irr}(n\!\!\:\cdot\!\!\:
  A_1,k\!\!\:\cdot\!\!\: A_2)$}
is T-smooth or empty.
\end{corollary}

\begin{corollary}  With the assumptions of Corollary
  \ref{cor:Theorem 1}, let $S_1,\dots,S_r$ be ordinary singularities
  of order $m_1,\dots,m_r$. Then \mbox{$V_{|dD_0|}^\text{\it
      irr}(S_1,\dots,S_r)$} is T-smooth or empty if
\begin{equation}
4\cdot \#(\mbox{nodes\/}) +
\sum_{m_i\geq 3}  \left\lfloor
    \frac{m_i^2\!+\!\!\:2m_i\!+\!\!\:5}{4}\right\rfloor^2
\! < \bigl(D_0^2- A(\Sigma,D_0)\bigr) d^2 - 2(D_0.K_\Sigma)  d \,.
\label{en12}\end{equation}
\end{corollary}

\medskip
\noindent
In our forthcoming book \cite{GLSBook2}, special emphasis is put on two
important examples, ESF of curves on a smooth hypersurface
\mbox{$\Sigma \subset \P^3$}:

\begin{theorem}[\cite{GLSBook2}]
\label{Theorem 4}
Let \mbox{$\Sigma \subset \P^3$} be a smooth hypersurface
of degree \mbox{$e\geq 5$}, let \mbox{$D_0=(e-4)H$}, $H$ a
hyperplane section, and let
\mbox{$d\in \Q$}, \mbox{$d > 1 $}. Suppose, moreover, that
\mbox{$\max_i \tci'(S_i)\leq d\!\!\:-\!\!\:2$}. If
\begin{equation}
\label{Nr neu}
 \sum_{i=1}^r \tci'(S_i) \,<\,
 \frac{e(e\!\!\:-\!\!\:4)^2}{4}\cdot  (d\!\!\:-\!\!\:1)^2
\end{equation}
and
\begin{equation}
 \label{cond 2}
\sum_{i=1}^r \frac{\bigl(\tci'(S_i)\!\!\:
  +\!\!\:1\bigr)^2}{1\!\!\:-\!\!\:\frac{\tci'(S_i)+1}{d}} \, < \,
  e(e\!\!\:-\!\!\:4)^2 \cdot d^2
\end{equation}
then \mbox{$V_{|dD_0|}^{\text{\it irr}}(S_1,\dots,S_r)$} is T-smooth
or empty.
\end{theorem}

\noindent
For ESF of nodal curves, respectively for ESF of curves with only
nodes and cusps, we can 
easily conclude the following quadratic sufficient conditions for
T-smoothness:

\begin{corollary}\label{cn2}
Let \mbox{$\Sigma \subset \P^3$} be a smooth hypersurface
of degree \mbox{$e\geq 5$}, let {$D_0=(e-4)H$}, $H$ a
hyperplane section, and let
\mbox{$d\in \Q$}.
\begin{enumerate}
\item[(a)] If   \mbox{$d \geq 3$} and
\begin{equation}
\label{yyy nodes}
 4\!\:n \:< \: e(e\!\!\:-\!\!\:4)^2\cdot d(d\!\!\:-\!\!\:2)
\end{equation}
then \mbox{$V_{|dD_0|}^{\text{\it irr}} (n\!\!\:\cdot\!\!\: A_1) $}
is T-smooth or empty.
\item[(b)] If \mbox{$d \geq 4$} and
\begin{equation}
\label{yyy}
 4\!\:n+ 9\!\:  k
\:< \: e(e\!\!\:-\!\!\:4)^2\cdot d(d\!\!\:-\!\!\:3)
\end{equation}
then \mbox{$V_{|dD_0|}^{\text{\it irr}} (n\!\!\:\cdot\!\!\: A_1, k\!\!\:\cdot
\!\!\:A_2) $}  is T-smooth or empty.
\end{enumerate}
\end{corollary}

\noindent
In the case of a quintic surface \mbox{$\Sigma\subset \P^3$}
 (\mbox{$e=5$}), Chiantini and Sernesi \cite{ChS}
provide examples of curves
\mbox{$C\subset \Sigma$},  \mbox{$C\equiv d\!\!\:\cdot\!\!\:H$},
\mbox{$d\geq 6$},  having
\mbox{$(5/4)\!\!\:\cdot\!\!\:(d\!\!\:-\!\!\:1)^2$} nodes such that
 \mbox{$V_{|C|}^{\text{\it irr}} (n\!\!\:\cdot\!\!\: A_1) $}
is not T-smooth at $C$. In particular, these examples show
that the exponent $2$ for $d$ in
the right-hand side of (\ref{yyy nodes}) is the best possible. Actually, it
even shows that for families of nodal curves on a quintic surface the condition
 (\ref{yyy nodes}) is asymptotically exact.

Performing more thorough computations in the proof of
\cite[Thm.\ 1]{GLS97}, Keilen improved the result of Theorem
\ref{Theorem 1} for surfaces with Picard number one or two. For
example, the following statement generalizes Theorem \ref{Theorem
4}:

\begin{theorem}[{\cite{Kei2}}]\label{tn1}
Let $\Sigma$ be a surface with Neron-Severi group
\mbox{$L\cdot\Z$}, $L$ being ample, let \mbox{$D=d\cdot L$}, let
$S_1,\dots,S_r$ be topological or analytic singularity types, and
let \mbox{$K_\Sigma=k_{\Sigma}\cdot L$}. Suppose that
\mbox{$d\ge\max\{k_{\Sigma}+1,-k_{\Sigma}\}$}, and
\begin{equation}\sum_{i=1}^r\gamma'(S_i)<\alpha\cdot(D-K_\Sigma)^2,\quad
\alpha=\frac{1}{\max\{1,1+k_{\Sigma}\}}\,.\label{en6}\end{equation} Then
either \mbox{$V^{\text{\it irr}}_{\Sigma,|D|}(S_1,\dots,S_r)$} is
empty or it is T-smooth.
\end{theorem}

\medskip\noindent
It is interesting that the same invariant $\gamma'(S)$
comes out of the proof using the Bogomolov-Reider theory of unstable
rank two vector bundles on surfaces instead of the Castelnuovo function
theory (as done in the proof of \mbox{Theorem \ref{Smoothness condition}}).

\subsection{ESF of Hypersurfaces}\label{sec:linear}

\noindent
{\bf Sufficient Conditions for T-Smoothness.}
Denote by \mbox{$V^n_d(S_1,\dots,S_r)$}
the set of hypersurfaces of degree
$d$ in $\P^n$, \mbox{$n\ge 3$}, whose singular locus consists of $r$
isolated singularities of analytic
types $S_1,\dots,S_r$, respectively. The following theorem was proved
independently by Shustin and Tyomkin and by Du Plessis and
Wall (using different methods of proof):

\begin{theorem}[\cite{ST,DPW}]\label{t7}
Let $S_1,\dots,S_r$ be analytic singularity types satisfying
\begin{equation}
\sum_{i=1}^r\tau(S_i)\,<\:\begin{cases}4d-4\;&\text{if}\ d\ge 5,\\
18\;&\text{if}\ d=4,\\
16\;&\text{if}\ d=3\ .\end{cases}\label{e42}
\end{equation}
Then the variety $V^n_d(S_1,\dots,S_r)$ is T-smooth, that is,
  smooth of the expected codimension $\sum_{i=1}^r\tau(S_i)$ in
  \mbox{$\big|H^0\bigl({\ko}_{\P^n}(d)\bigr)\big|$}.
\end{theorem}
A similar statement for topological types of singularities cannot be
true in general, because, for some singularities, the \mbox{$\mu=\const$}
stratum in a versal deformation base
(being a local topological ESF) is not smooth
\cite{Luengo1}. However, for semiquasihomogeneous singularities
the \mbox{$\mu=\const$} stratum
in a versal deformation base is smooth \cite{Var82}. 

Here a hypersurface singularity \mbox{$(W,z)\subset (\P^n\!,z)$} is
called {\em semiquasihomogeneous\/}\index{semiquasihomogeneous hypersurface singularity} 
({\em SQH\/}) if there are local analytic coordinates such that
$(W,z)$ is given by a power series 
\begin{equation}
f=\sum_{\wdeg(\balpha)\ge a} c_{\balpha}
  \bx^{\balpha}\in \C\{\bx\}\,,
\end{equation}
(where $\wdeg(\balpha):=\sum_{j=1}^n w_j\alpha_j$, \mbox{$\bw\in
  (\Z_{>0})^n$}) such that  
the Newton polygon of $f$ is convenient (that is, it intersects all  
coordinate axes), and such that the {\em principal part of $f$}, 
$$f_{0}= \sum_{\wdeg(\balpha)=a}c_{\balpha} \bx^{\balpha}\,,$$ defines an
isolated singularity at the origin. We introduce the ideal
$$
\Isqh(f) = \langle \bx^{\balpha} \mid \wdeg(\bx^{\balpha})\geq \wdeg(f_0)\rangle +
\Iea(f) \subset \C\{\bx\}\,.
$$
Its codimension in $\C\{\bx\}$ is an invariant of the topological type
$S$ of $(W,z)$. We denote it by $\deg \Zsqh(S)$. 

Using the smoothness of the \mbox{$\mu=\const$} stratum for SQH
singularities, in \cite{GLSBook2} we prove the following extension of
Theorem \ref{t7}:

\begin{theorem}[\cite{GLSBook2,ST}]\label{t8} 
Let $S_1,\dots,S_q$, $q<r$, be analytic singularity types, and let
$S_{q+1},\dots,S_r$ be topological types of semiquasihomogeneous
singularities. If
\begin{equation}\sum_{i=1}^q\tau(S_i)+\sum_{i=q+1}^r\deg\Zsqh(S_i)
<\, \begin{cases}4d-4\;&\text{if}\ d\ge 5,\\
18\;&\text{if}\ d=4,\\
16\;&\text{if}\ d=3\ ,\end{cases}\label{en10}\end{equation} then
the variety $V^n_d(S_1,\dots,S_r)$ is T-smooth.
\end{theorem}

\noindent
Du Plessis and Wall \cite{DPW} consider also
linear systems of hypersurfaces having a fixed intersection
with a hyperplane, which does not pass through the singular points, and
obtain the statement of Theorem \ref{t7} for such linear systems
under the condition
$$\sum_{i=1}^r\tau(S_i)<\, \begin{cases}3d-3\;&\text{if}\ d\ge 4,\\
8\;&\text{if}\ d=3\ .\end{cases}$$
Similarly, one can formulate an analogue of Theorem \ref{t8}.

\bigskip\noindent
{\bf  Non-T-Smooth ESF of Hypersurfaces.}
The (linear) condition (\ref{e42}) in Theorem \ref{t7} is not
necessary, as already seen in the case \mbox{$n=2$}.
In the following, we discuss which kind of sufficient conditions one might
expect.

It would be natural to extend or generalize the corresponding
results for plane curves  to higher
dimensions. Theorem \ref{t7} is a generalization of the $4d$-condition
for plane curves (see, e.g., \cite{GrL}).  The classical
$3d$-condition \eqref{k<3d}, however, cannot be
extended to higher dimensions in the same form. This follows, since
for plane curves it allows any number of nodes, while
there are surfaces of degree \mbox{$d\to\infty$} in $\P^3$ with
\mbox{$5d^3/12+O(d^2)$} nodes \cite{Chm} as only singularities. For
\mbox{$d\gg 1$}, these nodes
must be dependent as
{$\dim\big|H^0\bigl({\ko}_{\P^3}(d)\bigr)\big|=d^3/6+O(d^2)$}.
We also point out another important difference
between the case of curves and the case of higher dimensional
hypersurfaces. The quadratic numerical sufficient conditions for
T-smoothness of ESF of plane curves are close
to necessary conditions for the existence, which are quadratic in the degree
$d$ as well.
In higher dimensions the situation is different. Namely,
necessary conditions for the non-emptiness of $V_d^n(S_1,\dots,S_r)$,
such as
$$\sum_{i=1}^r\mu(S_i)\le(d-1)^n\ ,$$
are of order $d^n$ in the right-hand side and, for a fixed $n$,
there exist hypersurfaces with number of arbitrary singularities of
order $d^n$ (see Section \ref{sec:existESF of Hypersurfaces} below). 
However, any possible sufficient condition for
T-smoothness, in the form of an upper bound to the sum of certain
positive singularity invariants, can have at most a quadratic function
in $d$ on the right-hand side. Indeed, the following lemma
allows us to extend
examples of analytic ESF of plane
curves which are non-T-smooth to higher dimensions
such that the degree of the hypersurface and the total Tjurina number
are not changed.

\begin{lemma}[{\cite{ST}}]
\label{l8}
Let $C$ be a reduced plane curve of degree \mbox{$d>2$}.
Then, for any \mbox{$n>2$}, there exists a hypersurface
\mbox{$W\subset \P^n$} of
degree $d$ having only isolated singular points, such that
\mbox{$\tau(W)=\tau(C)$} and
$$h^1\bigl({\kj}_{\Zea(C)/\P^2}(d)\bigr)=
h^1\bigl({\kj}_{\Zea(W)/\P^n}(d)\bigr)\,.$$
\end{lemma}

\noindent
Notice that, in view of the obstructed families given in \cite{DPW},
this also yields that the inequality (\ref{e42}) cannot be
improved by adding a constant.

\section{Independence of Simultaneous Deformations}\label{sec:indep}

The T-smoothness problem for equisingular families is closely
related to the independence of simultaneous deformations of
isolated singular points of a curve on a surface, or of a
hypersurface in a smooth projective algebraic variety, and we present here
sufficient conditions for the independence of simultaneous
deformations which are analogous to the T-smoothness criteria of Section
\ref{sec:smoothness}.

\subsection{Joint Versal Deformations}
Let $W$ be a hypersurface with $r$ isolated singularities of analytic
types $S_1,\dots,S_r$, lying in a smooth projective algebraic variety
$X$. The obstructions to the versality of 
the joint deformation of the singularities of $W$ induced by
the linear system $|W|$ lie in the group $H^1({\cal
J}_{\Zea(W)/X}(W))$. In turn, the latter group is the obstruction to the
T-smoothness of the germ at $W$ of the ESF of all hypersurfaces
\mbox{$W'\in |W|$} having precisely $r$ singularities of analytic
types $S_1,\dots,S_r$ (see Proposition \ref{H0H1prop} for the
case that $X$ is a surface).

Thus, the aforementioned sufficient conditions for T-smoothness,
formulated in Theorems \ref{Smoothness Theorem}, \ref{Theorem 1},
\ref{Theorem 4}, \ref{tn1}, \ref{t7} and in
Corollaries \ref{Corollary 3.4}, \ref{cor:Theorem 1}, \ref{cn1},
\ref{cn2}, are also sufficient conditions for the versality of the
joint deformations of the singular points:

\begin{theorem}\label{thm:indep}
Let $W$ be a hypersurface with only isolated singular points
$z_1,\dots,z_r$ in a smooth projective algebraic variety $X$ of dimension
\mbox{$n\ge 2$}. Then the germ at $W$ of the linear system $|W|$ induces a
joint versal deformation of all the singularities of $W$, if one
of the following conditions holds:
\begin{enumerate}\item[(i)] \mbox{$X=\P^2$}, and \mbox{$W=C$} is an
  irreducible
  curve of degree $d$ such that (\ref{Smoothness condition}) is
  satisfied with \mbox{$\gamma'=\gea$}; if $C$ has
  $n$ nodes and $k$ cusps as only singularities, condition
  (\ref{Smoothness condition}) can be replaced by (\ref{en1});
\item[(ii)]
\mbox{$X=\Sigma$} is a surface with Neron-Severi
  group   \mbox{$L\cdot\Z$}, 
$L$ being ample, and with canonical divisor \mbox{$K_\Sigma=k_\Sigma\cdot L$};
\mbox{$W=C\sim d\cdot L$} is an 
irreducible curve, satisfying \mbox{$d\ge\max\{k_\Sigma+1,-k_\Sigma\}$} and
condition (\ref{en6}) with \mbox{$\gamma'=\gea$}; 
\item[(iii)] \mbox{$X=\Sigma$} is a
surface, $W=C$ is an irreducible curve such that $C$,
\mbox{$C\!\!\:-\!\!\:K_{\Sigma}$} are ample and \mbox{$C^2\geq
\max \,\{K_{\Sigma}^2, A(\Sigma,C)\}$}, and $\Sigma, C$ satisfy
conditions (\ref{smoot1.0}) and (\ref{smoot1.1}) with \mbox{$\tci'=\teaci$};
\item[(iv)] \mbox{$X=\Sigma$} is a surface with $-K_\Sigma$ nef, \mbox{$W=C$}
  is an irreducible curve such that \mbox{$C\sim dD_0$}, $D_0$ an ample
divisor on $\Sigma$ such that \mbox{$D_0^2\geq A(\Sigma,D_0)$},
\mbox{$d>0$} satisfying $d^2\geq K_\Sigma^2/D_0^2$ and condition
(\ref{smoot1.1'}) with \mbox{$\tci'=\teaci$} (the
latter condition reduces to (\ref{smoot1.1'nc}) if $C$ has $n$
nodes and $k$ cusps as only singularities);
\item[(v)] \mbox{$X=\Sigma \subset \P^3$} is a hypersurface of degree
  \mbox{$e\geq 5$}, \mbox{$W=C$} is an irreducible curve such that \mbox{$C\sim
    d(e-4)H$}, where $H$ is a hyperplane section, \mbox{$d\in \Q$},
  \mbox{$d>1$}, such 
 that \mbox{$\max_i \teaci(C,z_i)\leq d\!\!\:-\!\!\:2$}, and such that the
 conditions (\ref{Nr neu}), (\ref{cond 2}) are satisfied with
 \mbox{$\tci'=\teaci$} (the latter two conditions turn into
 (\ref{yyy}), if $C$ has $n$ nodes and $k$ cusps as only
 singularities); 
\item[(vi)] \mbox{$X=\P^n$}, \mbox{$n\ge 2$}, $W$ is a reduced hypersurface of
degree $d$, and condition (\ref{e42}) is fulfilled.
\end{enumerate}
\end{theorem}

\subsection{Independence of Lower Deformations}

Also the T-smoothness of a topological equisingular family has a
deformation theoretic counterpart: the independence of lower deformations of
isolated singularities.

Let $W$ be a hypersurface with only isolated singular points
$z_1,\dots,z_r$ in a smooth projective algebraic variety $X$ of dimension
\mbox{$n\ge 2$}. For sake of simplicity, we assume that the singular points
of $W$ are all semiquasihomogeneous (SQH). That is, for each singular
point $z_i$, there are local analytic 
coordinates \mbox{$\bx=(x_1,\dots,x_n)$} on $X$ such that the germ $(W,z_i)$ is
given by 
\begin{equation}F_i=\sum_{\ell_i(\balpha)\ge a_{i}}A^{(i)}_{\balpha}
  \bx^{\balpha}\in \C\{\bx\}=\C\{x_1,\dots,x_n\}\, 
,\label{en8}\end{equation} 
where \mbox{$\ell_{i}(\balpha)= \sum_{j=1}^n
  w_j^{(i)}\alpha_j$} for some
\mbox{$\bw^{(i)}\!=(w^{(i)}_1\!,\dots,w^{(i)}_n)\in (\Z_{> 0})^n$}, 
 the Newton polygon 
  of $F_i$ is convenient and \mbox{$F_{i,0}=
\sum_{\ell_i(\balpha)=a_i}A^{(i)}_{\balpha} \bx^{\balpha}$} defines an
isolated singularity at the origin.
We call $F_i$ a {\em SQH representative\/} of $(W,z_i)$ with {\em principal
part\/} $F_{i,0}$.
Note that the class of SQH singularities includes all simple
singularities. 

We fix SQH representatives \mbox{$F_1,\dots,F_r$} for 
$(W,z_1),\dots,(W,z_r)$. Then by a {\em deformation
  pattern\/} for $(W,z_i)$ (respectively for
$F_i$), we denote any affine hypersurface of $\C^n$ given by a polynomial
$$G_i= \sum_{0\leq \ell_i(\balpha)\le a_i} A^{(i)}_{\balpha}
\bx^{\balpha} = F_{i,0} + \sum_{0\leq \ell_i(\balpha)<
  a_i} A^{(i)}_{\balpha} 
\bx^{\balpha}  $$
and having only isolated singularities.
Note that the set of all deformation patterns for $(W,z_i)$ can be
identified with the subset 
\begin{equation}
\mathcal{P}_i= \left\{ (A^{(i)}_{\balpha})_{0\leq \ell_i(\balpha)<
    a_i} \,\left|\: \begin{array}{c} A^{(i)}_{\balpha} \!\in \C \:\text{ and 
    the hypersurface } \{G_i=0\} \\\text{ has only isolated
    singularities in } 
    \C^n
  \end{array}
\right.\right\}\label{en9} 
\end{equation} 
of the affine space parametrized by all lower coefficients. 

After applying a weighted homothety \mbox{$(x_1,\dots,x_n) \mapsto
  (\lambda^{w^{(i)}_{1}}x_1,\dots,\lambda^{w^{(i)}_{n}}x_n)$} (with \mbox{$|\lambda|$}
sufficiently large) to a given deformation pattern, we
can assume that the coefficients $A^{(i)}_{\balpha}$,
\mbox{$\ell_i(\balpha)<a_i$}, are sufficiently small such that, for a
fixed large closed ball \mbox{$B_i\subset\C^n$}, the intersection 
\mbox{$\{G_i=0\}\cap \partial B_i$} is close to \mbox{$\{F_{i,0}=0\}\cap 
\partial B_i$}. For each \mbox{$i=1,\dots,r$}, we choose a small
closed regular neighbourhood 
$V_i$ of $z_i$ in $X$ and a $C^\infty$-diffeomorphism
$$\varphi_i:(\partial B_i,\{F_{i,0}=0\}\cap \partial B_i)\to
(\partial V_i,W\cap\partial V_i)$$ which is close to a weighted
homothety in the coordinates $\bx$.

If $p_1,\dots,p_s$ are the singular points of
\mbox{$\{G_i=0\}$} in $B_i$ then, for each $j$, we choose the
topological\footnote{When talking about topological types, we always
  assume \mbox{$n=2$}.} or the analytic equivalence relation, and
by the germ of the {\em equisingularity stratum in $\kp_i$
at\/ $G_i$}\index{equisingularity stratum} we mean the germ of the set of all deformation patterns
$G'_i$ having singular points  $p'_1,\dots,p'_s$ close to
$p_1,\dots,p_s$ such that $(G_i,p_i)$ and $(G'_i,p'_i)$ have the same
type (with respect to the chosen equivalence).
We call the deformation pattern defined by $G_i$ {\em transversal\/}
if the germ at $G_i$ of the equisingularity stratum in the space ${\cal
P}_i$ is T-smooth.  

Given a one-parameter deformation \mbox{$W_t\in|W|$}, \mbox{$t\in(\C,0)$}, of
the hypersurface \mbox{$W=W_0$} such that \mbox{$W_t\cap V_i$},
\mbox{$t\ne 0$}, is equisingular (with respect to the fixed
equivalences), we say that $\{W_t\}$ {\em matches the deformation
  pattern\/}\index{deformation pattern}
$G_i$ for $(W,z_i)$ if, for \mbox{$t\ne 0$}, there is a homeomorphism
$\psi_{i,t}$ of \mbox{$(B_i,\{G_i=0\}\cap B_i)$} onto
\mbox{$(V_i,W_t\cap V_i)$}. Moreover, at
those singular points $p_j$ of \mbox{$\{G_i=0\}$} where the analytic
equivalence was chosen, we additionally require that $\psi_{i,t}$ induces
an analytic isomorphism in a neighbourhood  of $p_j$.

\medskip
\noindent
The main result of \cite{Sh7} is

\begin{theorem}\label{tn10}
Let $W$ be a hypersurface in a smooth projective algebraic variety $X$ of
dimension $n\ge 2$ with only isolated semiquasihomogeneous singular
points $z_1,\dots,z_r$. If the germ at $W$ of the topological ESF in the
linear system $|W|$ is T-smooth, then, for each tuple of independently
prescribed transversal deformation patterns 
for $(W,z_1),\dots,(W,z_r)$, there exists a one-parameter deformation
\mbox{$W_t\in |W|$}, \mbox{$t\in (\C,0)$}, of $W$ which matches the given
patterns.

Moreover, if all the given data are defined over the reals, then the
deformation $W_t$ and the matching homeomorphisms
$\psi_{i,t}$ can be chosen over the reals, too.
\end{theorem}

\noindent
The proof is constructive and is based on a version of the
patchworking method\index{patchworking method}, which is one of the main tools used
for finding 
sufficient conditions for the existence of hypersurfaces with
prescribed singularities (see Section \ref{sec:existence} below).

Various T-smoothness criteria for topological ESF immediately
imply numerical sufficient conditions for the independence of
one-parameter deformations matching given deformation patterns.
Skipping the linear conditions given in
\cite{GrK,GrL,Shu87,Shu91,Vas}, we collect in the following theorem
the outcome of the criteria in Section \ref{sec:smoothness}. Notice that
the case that $W$ has only simple singularities is already covered by
Theorem \ref{thm:indep}. 

\begin{theorem} Let $W$ be a hypersurface in a smooth projective
  algebraic variety $X$ of dimension $n\ge 2$ with only isolated
  semiquasihomogeneous singular points $z_1,\dots,z_r$. Moreover, let
  $G_1,\dots,G_r$ define transversal deformation patterns for
  $(W,z_1),\dots,(W,z_r)$.
Then there exists a one-parameter deformation \mbox{$W_t\in |W|$},
\mbox{$t\in (\C,0)$}, of $W$ which matches these deformation patterns,
if one of the following conditions is satisfied:
\begin{enumerate}
\item[(i)] \mbox{$X=\P^2$} and \mbox{$W=C$} is an irreducible
  curve of degree $d$ such that (\ref{Smoothness condition}) holds
  with \mbox{$\gamma'=\ges$}; if
  $C$ has only ordinary singularities of order 
$m_1,\dots,m_r$,
  condition (\ref{Smoothness condition}) can be replaced by (\ref{en11});
\item[(ii)] \mbox{$X=\Sigma$} is a surface with Neron-Severi
  group   \mbox{$L\cdot\Z$}, 
$L$ being ample, and with canonical divisor \mbox{$K_\Sigma=k_\Sigma\cdot L$};
\mbox{$W=C\sim d\cdot L$} is an 
irreducible curve, satisfying \mbox{$d\ge\max\{k_\Sigma+1,-k_\Sigma\}$} and
condition (\ref{en6}) with \mbox{$\gamma'=\ges$};
\item[(iii)]  \mbox{$X=\Sigma$} is a
surface, $W=C$ is an irreducible curve such that $C$,
\mbox{$C\!\!\:-\!\!\:K_{\Sigma}$} are ample and \mbox{$C^2\geq
\max \,\{K_{\Sigma}^2, A(\Sigma,C)\}$}, and $\Sigma, C$ satisfy
conditions (\ref{smoot1.0}) and (\ref{smoot1.1}) with 
\mbox{$\tci'=\tesci$} (conditions (\ref{smoot1.0}),
(\ref{smoot1.1}) can be replaced by (\ref{en12}) if $C$ has only
ordinary singularities of order $m_1,\dots,m_r$);
\item[(iv)]  \mbox{$X=\Sigma$} is a surface with $-K_\Sigma$ nef, \mbox{$W=C$}
  is an irreducible curve such that \mbox{$C\sim dD_0$}, $D_0$ an ample
divisor on $\Sigma$ with \mbox{$D_0^2\geq A(\Sigma,D_0)$},
\mbox{$d>0$} satisfying \mbox{$d^2\geq K_\Sigma^2/D_0^2$} and condition
(\ref{smoot1.1'}) with \mbox{$\tci'=\tesci$};
\item[(v)]  \mbox{$X=\Sigma \subset \P^3$} is a hypersurface of degree
  \mbox{$e\geq 5$}, \mbox{$W=C$} is an irreducible curve such that \mbox{$C\sim
    d(e-4)H$}, where $H$ is a hyperplane section, \mbox{$d\in \Q$},
  \mbox{$d>1$}, such 
 that \mbox{$\max_i \tesci(C,z_i)\leq d\!\!\:-\!\!\:2$}, and such that the
 conditions (\ref{Nr neu}), (\ref{cond 2}) are satisfied with
 \mbox{$\tci'=\tesci$};
\item[(vi)]  \mbox{$X=\P^n$}, \mbox{$n\ge 2$}, $W$ is a reduced hypersurface of
degree $d$, and condition (\ref{en10}) is fulfilled with \mbox{$q=0$}.
\end{enumerate}
\end{theorem}

\section{Existence}  \label{sec:existence}

\noindent
In the following, we describe two methods which lead to general numerical sufficient
conditions for the existence of projective hypersurfaces with
prescribed singularities. Both approaches are based on the reduction
of the existence problem to an $H^1$-vanishing problem for the ideal sheaves of certain
zero-dimen\-sional schemes associated with topological, respectively analytic,
types of singularities. One way is to associate directly
a zero-dimensional scheme corresponding to the $r$ prescribed
singularity types $S_1,\dots,S_r$ (fixing the position of the
singular points) and to produce a sufficient condition by using an appropriate
$H^1$-vanishing criterion. Another way is to 
construct, first, a projective hypersurface with ordinary
singularities (in general position) and
then to deform it into a hypersurface with the prescribed
singularities using the patchworking construction
(\cite{Shu98,Sh13}). The sufficient 
conditions obtained by each of the two approaches do not cover the
conditions obtained by the other approach in general. Hence, both methods
are needed.

\subsection{ESF of Plane Curves}

\noindent
With a reduced plane curve germ \mbox{$(C,z)\subset (\P^2\!,z)$} we
associate the following zero-di\-men\-sional schemes of $\P^2$ (with
support $\{z\}$):
\begin{itemize}
\itemsep3pt
\item $Z^s(C,z)$, the {\em singularity scheme}, defined by the ideal
\begin{eqnarray*}
I^s(C,z) &:=&
\bigl\{g \in \ko_{\P^2,z} \:\big|\: \mt \widehat{g}_{(q)}
\geq \mt \widehat{C}_{(q)}\!\;\text{ for each } q \in \kt^\ast(C,z)\bigr\}\,,
\end{eqnarray*}
where $\kt^\ast(C,z)$ denotes the tree of essential infinitely near
points, and $\widehat{g}_{(q)}$ (resp.\ $\widehat{C}_{(q)}$) is the total
transform of $g$ (resp.\ of $(C,z)$) at $q$ 
(see \cite{GLS96} or \cite{GLSBook1} for details);
\item \mbox{$\Zsst(C,z):=Z^s(CL,z)$}, where $L$ is a curve which is
  smooth and transversal to (the tangent cone of) $C$ at $z$;
\item $Z^a(C,z)$, the scheme defined by the ideal
\mbox{$ I^a(C,z)\subset\ko_{\P^2,z}$} encoding the analytic type (see
\cite{GLS98} for a definition);
\item $\Zast(C,z)$, the scheme defined by the ideal
$\fm_z I^a(C,z)$, where $\fm_z$ denotes the maximal ideal of $\ko_{\P^2,z}$.
\end{itemize}

\bigskip\noindent
{\bf $\bH^{\bone}$-Vanishing Approach.} 
The following proposition allows us to deduce the existence of plane
curves with prescribed singularities from an $H^1$-vanishing statement:

\begin{proposition}[\cite{Sh01}]
\label{Existenz96}
(1) Given a zero-dimensional scheme
  \mbox{$Z\subset\P^2$}, a point \mbox{$z\in \P^2$} outside the support of $Z$
  and a reduced curve germ \mbox{$(C,z)\subset (\P^2\!,z)$} satisfying
\begin{equation}
H^1\bigl(\kj_{Z\cup \Zsst(C,z)/\P^2}(d)\bigr)=0\,.\label{ee1}
\end{equation}
Then there exists a curve \mbox{$D\in\big|H^0\bigl(\kj_{Z\cup
    Z^s(C,z)/\P^2}(d)\bigr)\big|$} such
that the germ of $D$ at $z$ is topologically equivalent to $(C,z)$.
Moreover, these curves $D$ form a dense open subset in
\mbox{$\big|H^0\bigl(\kj_{Z\cup Z^s(C,z)/\P^2}(d)\bigr)\big|$}.

\medskip\noindent
(2) In the previous notation, let
\begin{equation}
H^1\bigl(\kj_{Z\cup\Zast(C,z)/\P^2}(d)\bigr)=0\, .\label{ee2}
\end{equation}
Then there exists a curve \mbox{$D\in\big|H^0\bigl(\kj_{Z\cup
    Z^a(C,z)/\P^2}(d)\bigr)\big|$} such that the germ of $D$ at $z$ 
is analytically equivalent to $(C,z)$. These curves
$D$ form a dense open subset in \mbox{$\big|H^0\bigl(\kj_{Z\cup
    Z^a(C,z)/\P^2}(d)\bigr)\big|$}.
\end{proposition}

\noindent
Together with the $H^1$-vanishing theorem for generic zero-dimensional
  schemes given in \cite{Sh01} (using a Castelnuovo function
  approach), Proposition \ref{Existenz96} yields:

\begin{theorem}[\cite{Sh01}]
\label{te2}
Let $(C_1,z_1),\dots,(C_r,z_r)$ be reduced plane curve germs,
let $n$ be the number of nodes, $k$ the number of cusps and $t$ the number
of $A_{2m}$ singularities, \mbox{$m\ge 2$}, among
the singularities $(C_i,z_i)$, \mbox{$i=1,\dots,r$}.

\medskip\noindent
(1) If
\begin{equation}
6n+10k+\frac{49}{6}t+\frac{625}{48}
 \sum_{(C_i,z_i)\ne A_1,A_2}\!\!\del(C_i,z_i)\,\le\, d^2\!-2d+3\,
 ,\label{ee5}
 \end{equation}
then there exists a reduced, irreducible plane curve of degree $d$ having $r$
singular points topologically equivalent to $(C_1,z_1), \dots,
(C_r,z_r)$, respectively, as its only singularities.

\medskip\noindent
(2) If
 \begin{equation}
 6n+10k+\sum_{(C_i,z_i)\ne
 A_1,A_2}\!\!\frac{\bigl(7\mu(C_i,z_i)+2\del(C_i,z_i)\bigr)^2}{6
\mu(C_i,z_i)+3\del(C_i,z_i)}\,\le\, d^2\!-2d+3\, ,\label{ee6}
 \end{equation}
then there exists a reduced, irreducible plane curve of degree $d$ having $r$
singular points analytically equivalent to
$(C_1,z_1), \dots ,(C_r,z_r)$, respectively, as its only
singularities.
\end{theorem}

\noindent
See \cite{Sh01} for a slightly stronger result. Note that
condition (\ref{ee6}) can be weakened to the following simple
form (in view of \mbox{$\del\le 3\mu/4$} for reduced plane curve
singularities different from nodes):
\begin{equation}
\label{eqn:d^2/9}
\sum_{i=1}^r\mu(C_i,z_i)\le\frac{1}{9}(d^2\!-2d+3)\, .
\end{equation}

\noindent
Next, we pay special attention to the case of curves with exactly one
singular point, because such curves are an essential ingredient for the
patchworking approach to the existence problem.

\medskip\noindent
{\bf Curves with one Singular Point and Order of T-existence.} 
Let \mbox{$(C,z)$} be a reduced plane curve
singularity.  Denote by
$e^s(C,z)$, resp.\ $e^a(C,z)$,  the minimal degree $d$ of a
plane curve \mbox{$D\subset \P^2$} whose singular locus consists of a unique
point $w$ such that $(D,w)$ is topologically (resp.\ analytically)
equivalent to $(C,z)$ and which satisfies the condition
\begin{equation}
H^1\bigl(\kj_{\Zes(D)/\P^2}(d-1)\bigr)=0\,, \: \text{ resp.\ }\,
H^1\bigl(\kj_{\Zea(D)/\P^2}(d-1)\bigr)=0\,. \label{ee20}
\end{equation}
We call $e^s(C,z)$ (resp.\ $e^a(C,z)$) the {\em order of
  T-existence\/}\index{order of T-existence} for the topological (resp.\
analytic) singularity type represented by $(C,z)$.

\begin{lemma}\label{le13}
(1) Let $D$ be a plane curve as in the definition of the order of T-existence, and
  let $L$ be a straight line which does not pass through the singular
  point $w$ of $D$. Then the
  germ at $D$ of the family of curves of degree $d$ having in a
  neighbourhood of $w$ a singular point which is topologically (respectively
  analytically) equivalent to $(C,z)$ is smooth of the expected
  dimension, and it intersects transversally the linear system
$$\Bigl\{G\in\big|H^0\bigl(\ko_{\P^2}(d)\bigr)\big|\ \Big|\ G\cap
L=D\cap L\Bigr\}\ .$$

\smallskip\noindent
(2) Let \mbox{$L\subset\P^2$} be a straight line. Then
the set of $d$-tuples $(z_1,\dots,z_d)$ of distinct points on $L$ for which
there is a curve $D$ of degree $d$ as in the definition of
the order of T-existence satisfying \mbox{$D\cap L=\{z_1,\dots,z_d\}$}
is Zariski open in $\Sym^d(L)$.
\end{lemma}

\medskip\noindent
Combining Theorem \ref{te2} and the existence result for plane curves
with simple singularities in \cite{Lossen}, we get the following
estimates for $e^s$ and $e^a$:

\begin{theorem}\label{te3}
 If $(C,z)$ is a simple plane curve singularity then
$$ e^s(C,z)=e^a(C,z)
\left\{
\renewcommand{\arraystretch}{1.2}
  \begin{array}{ll}
  \le 2\lfloor\sqrt{\mu\!\!\:+\!\!\:5}\rfloor & \text{if $(C,z)$ of
    type $A_\mu$, $\mu\ge 1$}\,,\\
  \le 2\lfloor\sqrt{\mu\!\!\:+\!\!\:7}\rfloor\!\!\:+\!\!\:1 & \text{if
    $(C,z)$ of type $D_\mu$, $\mu\ge
    4$}\,,\\
  =\left\lfloor \mu/2\right\rfloor+1 & \text{if $(C,z)$  of
    type $E_\mu$, $\mu=6,7,8$}.
  \end{array}
\right.
$$
 If $(C,z)$ is not simple, then
 \begin{eqnarray*}
 e^s(C,z)& \le & \frac{25}{4\sqrt{3}}\sqrt{\del(C,z)}-1\, ,\\
 e^a(C,z)& \le &
 \frac{7\mu(C,z)+2\del(C,z)}{\sqrt{6\mu(C,z)+3\del(C,z)}}-1
\le 3\sqrt{\mu(C,z)}-1\, .
 \end{eqnarray*}
\end{theorem}

\noindent
For simple singularities with small Milnor number, these estimates are
far from being sharp. For instance, it is well-known that
$$e^s(A_\mu)= \left\{
\begin{array}{ll}
 4 & \text{ if $3\le\mu\le 7$}\,,\\
 5 & \text{ if $8\le\mu\le 13$}\,,\\
 6 & \text{ if $14\le\mu\le 19$}\,,
\end{array}
\right.\quad
e^s(D_\mu)= \left\{
\begin{array}{ll}
 3 & \text{ if $\mu=4$}\,,\\
 4 & \text{ if $\mu=5$}\,,\\
 5 & \text{ if $6\le\mu\le 10$}\,,\\
 6 & \text{ if $11\le\mu\le 13$}\,.\\
\end{array}
\right.
$$
Moreover, \mbox{$e^s(E_6)=e^s(E_7)=4$}, \mbox{$e^s(E_8)=5$}.

\bigskip\noindent
{\bf Curves with many Singular Points (Patchworking
  Approach).}\index{patchworking method} The following 
proposition is a special case of Proposition \ref{pe5} below, which is proved
by a reasoning based on patchworking:

\begin{proposition}\label{pe4}
Let $(C_1,z_1),\dots,(C_r,z_r)$ be reduced plane curve
singularities. Let \mbox{$m_i:=e^s(C_i,z_i)$} (resp.\
\mbox{$m_i:=e^a(C_i,z_i)$}), and assume that
\begin{equation}
H^1\bigl(\kj_{Z(\bm)/\P^2}(d-1)\bigr)=0\, ,\label{ee15}
\end{equation}
where $Z(\bm)$ is the fat point scheme supported at the
(generic) points \mbox{$p_1,\dots,p_r$} of $\P^2$ and defined by the
ideals $\fm_{p_i}^{m_i}$. Moreover, let
$$d>\max_{1\le i\le
r}e^s(C_i,z_i)\,,\quad\big(\text{resp.\ }\, d>\max_{1\le i\le
r}e^a(C_i,z_i)\,\big)\, .$$
Then there exists a reduced, irreducible plane
curve \mbox{$D\subset \P^2$} of degree $d$ with
\mbox{$\Sing(D)=\{p_1,\dots,p_r\}$} such that each germ $(D,p_i)$ is
topologically (resp.\ analytically) equivalent to $(C_i,z_i)$,
\mbox{$i=1,\dots,r$}.
\end{proposition}

\medskip\noindent
Applying the $H^1$-vanishing criterion of \cite[Thm. 3]{GengXu95}, we
immediately derive

\begin{corollary}\label{ce3}
Let $(C_1,z_1),\dots,(C_r,z_r)$ be reduced plane curve singularities
such that \mbox{$e^s(C_1,z_1)\ge\ldots\ge e^s(C_r,z_r)$}. If
\begin{eqnarray*}
e^s(C_1,z_1)+e^s(C_2,z_2)&\le&
d-1\,,\quad\ \text{as }\, r\ge 2\, ,\\
e^s(C_1,z_1)+\ldots+e^s(C_5,z_5)&\le& 2d-2\,,\quad\text{as }\, r\ge
5\,, \\
\sum_{i=1}^r(e^s(C_i,z_i)+1)^2 &<& \frac{9}{10}(d+2)^2\, ,
\end{eqnarray*}
then there exists a reduced, irreducible plane curve $C$ of degree $d$
with exactly $r$ singular points $p_1,\dots,p_r$, such that each germ
$(C,p_i)$ is topologically equivalent to $(C_i,z_i)$,
\mbox{$i=1,\dots,r$}.

The same statement holds true if we replace $e^s$ by $e^a$ and
the topological equivalence relation by the analytic one.
\end{corollary}

\noindent
Comparing the sufficient conditions of Theorem \ref{te2} and of Corollary
\ref{ce3}, we see that the existence criterion of Theorem
\ref{te2} is better for non-simple singularities,
whereas the criterion obtained from Corollary \ref{ce3} and the
estimates in Theorem \ref{te3} is better for simple
singularities.

\subsection{ESF of Curves on Smooth Projective Surfaces}

In \cite{KeT}, the following sufficient criterion for the
existence of curves with prescribed singularities is proved:

\begin{proposition}[\cite{KeT}]\label{pe5}
Let $\Sig$ be a smooth projective algebraic surface, $D$ a divisor
on $\Sig$, and \mbox{$L\subset\Sig$} a very ample divisor. Let
$(C_1,z_1), \dots, (C_r,z_r)$ be reduced plane curve singularities. Let
\begin{eqnarray}
&H^1\bigl(\kj_{Z(\bm)/\Sig}(D-L)\bigr)=0\ ,\label{ee21}\\
&\max\limits_{1\le i\le r}m_i<L(D-L-K_\Sig)-1\ ,\label{enn1}
\end{eqnarray}
where \mbox{$Z(\bm)\subset \Sigma$} is the fat point scheme supported at some (generic)
points \mbox{$p_1,\dots,p_r\in \Sigma$} and defined by the
ideals \mbox{$\fm_{p_i}^{m_i}$}, with \mbox{$m_i=e^s(C_i,z_i)$},
respectively \mbox{$m_i=e^a(C_i,z_i)$},  \mbox{$i=1,\dots,r$}.
Then there exists an irreducible curve \mbox{$C\in|D|$} such that
\mbox{$\Sing(C)=\{p_1,\dots,p_r\}$} and each germ $(C,p_i)$ is
topologically, resp.\ analytically, equivalent to
$(C_i,z_i)$, \mbox{$i=1,\dots,r$}.
\end{proposition}

\noindent
Combining this with the $H^1$-vanishing criterion of
\cite[Cor.\ 4.2]{KeT} and with the estimates for $e^s$, $e^a$ in
Theorem \ref{te3}, we get the following explicit
numerical existence criterion:

\begin{theorem}[\cite{GLSBook2}]
\label{te6}
Let $\Sig$ be a smooth  projective algebraic surface, $D$ a divisor
on $\Sig$ with \mbox{$D-K_{\Sig}$} nef, and \mbox{$L\subset\Sig$} a
very ample divisor. Let $(C_1,z_1), \dots, (C_r,z_r)$ be reduced plane
curve singularities, among them $n$ nodes and $k$ cusps.

\medskip\noindent
(1) If
\begin{eqnarray}
18n+32k+\frac{625}{24}\sum_{\del(C_i,z_i)>1}\!\!\del(C_i,z_i)&\,\le
\,&(D-K_\Sig-L)^2\, ,\label{ee24}\\
\frac{25}{4\sqrt{3}}\max_{1\le i\le r}\sqrt{\del(
C_i,z_i)}+1 &<&(D-L-K_\Sig).L\, ,\label{ee22}
\end{eqnarray}
 and, for each irreducible curve $B$ with \mbox{$B^2=0$} and
 \mbox{$\dim|B|_a>0$},
\begin{equation}
\frac{25}{4\sqrt{3}}\max_{1\le i\le r}\sqrt{\del(C_i,z_i)}\,<\,
 (D-K_\Sig-L).B+1\, ,\label{ee23}
\end{equation}
 then there exists a reduced, irreducible curve \mbox{$C\in|D|$} with
 $r$  singular points topologically equivalent to $(C_1,z_1), \dots,
 (C_r,z_r)$, respectively, as its only singularities.

\medskip\noindent
(2) If
\begin{eqnarray}
18n+32k+18\sum_{\mu(C_i,z_i)>2}\!\!\mu(C_i,z_i)&\,\le
\,&(D-K_\Sig-L)^2\,,\label{ee25}\\
3\max_{1\le i\le r}\sqrt{\mu(C_i,z_i)}+1&<&(D-L-K_\Sig).L\,
,\label{ee26}
\end{eqnarray}
 and, for each irreducible curve $B$ with \mbox{$B^2=0$} and
 \mbox{$\dim|B|_a>0$},
 \begin{equation}
3\max_{1\le i\le r}\sqrt{\mu(C_i,z_i)}\,<\, (D-K_\Sig-L).B+1\,
,\label{ee27}
\end{equation}
 then there exists a reduced, irreducible curve \mbox{$C\in|D|$} with
 $r$ singular points analytically equivalent to $(C_1,z_1), \dots,
 (C_r,z_r)$, respectively, as its only singularities.
\end{theorem}

\noindent Here $|B|_a$ means the family of curves algebraically
equivalent to $B$. A discussion of the hypotheses of Theorem
\ref{te6} for specific classes of surfaces as well as concrete
examples can be found in \cite{KeT}.

\subsection[ESF of Hypersurfaces in $\P^n$]{ESF of Hypersurfaces in
$\mathbf{P}^{\bn}$}\label{sec:existESF of Hypersurfaces} For
singular hypersurfaces in $\P^n$, \mbox{$n\ge 3$}, no general
asymptotically proper sufficient condition for the existence of
hypersurfaces with prescribed singularities (such as (\ref{eqn:d^2/9})
in the case of plane curves) is known. But, restricting ourselves
to the case of only simple singularities, in \cite{ShuWes04} even an asymptotically
optimal condition is given. 

To formulate this result, we need some notation: let ${\cal S}$ be a
finite set of analytic types of isolated hypersurface
singularities in $\P^n$. Define 
$$\alpha_n({\cal S})=\lim_{d\to\infty}\sup_{W_d}\frac{\tau(W_d)}{d^n}\, ,$$ 
where $W_d$ runs over the set of all hypersurfaces \mbox{$W_d\subset\P^n$} of
degree $d$ whose singularities are of types \mbox{$S\in{\cal S}$} and
which belong to the T-smooth component of the corresponding ESF. Here,
$\tau(W_d)$ stands for the sum of the Tjurina numbers $\tau(W_d,z)$
over all points \mbox{$z\in \Sing(W_d)$}. 
By $\alpha^\R_n({\cal S})$ we denote the respective
limit taken over hypersurfaces having only real singular points
of real singularity types \mbox{$S\in{\cal S}$}. Clearly,
\mbox{$\alpha^{\R}_n({\cal S})\le\alpha_n({\cal S})\le1/n!$}.

\begin{theorem}[\cite{ShuWes04,West,West1}]\label{thm:shuWes}
Let \mbox{$n\geq 2$}.
\begin{enumerate}
\item[(1)] For each finite set ${\cal S}$ of simple hypersurface
singularities in $\P^n$, 
we have $$\alpha^\R_n({\cal S})=\alpha_n({\cal
S})=\frac{1}{n!}\, .$$
\item[(2)] For each finite set ${\cal S}$ of analytic types of
  isolated hypersurface singularities of corank $2$ in $\P^n$, we have
$$\alpha_n({\cal S})\ge\alpha^\R_n({\cal
S})\ge\frac{1}{9n!}\, .$$
\end{enumerate}
\end{theorem}

\noindent
The proof exploits again the patchworking construction. It is
based on the following fact: for each simple singularity type $S$
and each \mbox{$n\ge 2$}, there exist an $n$-dimensional convex lattice
polytope $\Delta_n(S)$ of volume $\mu(S)/n!$ and a polynomial
\mbox{$F\in\C[x_1,\dots,x_n]$} (resp.\ \mbox{$F\in\R[x_1,\dots,x_n]$})
with Newton polytope $\Delta_n(S)$ which defines a hypersurface in the
toric variety $\Tor(\Delta_n(S))$ (associated to $\Delta_n(S)$) having
precisely one singular point of type $S$ in the torus $(\C^*)^n$
(resp.\ in the torus $(\R^*)^n$) and being
non-singular and transverse along the toric divisors in
$\Tor(\Delta_n(S))$.

\section{Irreducibility} \label{sec:irreducibility}

The question about the irreducibility of ESF
$\VC$ is more delicate than the existence and smoothness problem, in
particular, if one tries to find
sufficient conditions for the irreducibility. The results are by far
not that complete as for the other two problems. The irreducibility
problem is of special topological interest, since it is
connected with the problem of having within the same ESF different
fundamental groups of the complement of a plane algebraic curve.

As pointed out in the introduction, even the case of plane nodal
curves (Severi's conjecture) appeared to be very hard. The examples of
reducible ESF listed below indicate that for more complicated singularities,
beginning with cusps, possible numerical sufficient
conditions for the irreducibility should be rather different with
respect to their asymptotics to the
necessary existence conditions (as discussed in Section
\ref{sec:existence}).

\medskip \noindent
\mbox{\bf Approaches to the Irreducibility Problem.}
(1) One possible approach (for ESF of plane curves) consists of
  building for any two curves in the ESF a connecting
  path, using explicit equations of the curves, respectively
 of projective transformations. This method works for small degrees
only. Besides the classical case of conics and cubics, this method has
been used to prove that all ESF of quartic and quintic curves are
irreducible (cf.\ \cite{BrG,Wl2}). But for degrees \mbox{$d>5$}, this is no
longer true and the method is no more efficient (except for some very special
cases).

\medskip \noindent
(2) Arbarello and Cornalba \cite{ArC} suggested another approach. It
consists of relating the ESF to the moduli space of
plane curves of  given genus, which is known to be irreducible (cf.\
\cite{DeM}).  This gave some particular
results on families of plane nodal curves and plane curves with nodes and
cusps. Namely, Kang \cite{Kang} proved that the variety \Vdnk\  is
irreducible whenever
\begin{equation}
\label{Kang}
\frac{d^2\!-\!\!\:4d\!\!\:+\!\!\:1}{2} \,\leq \, n \,\leq \,
\frac{(d\!\!\:-\!\!\: 1)(d\!\!\:-\!\!\: 2)}{2}\,, \quad \;\; k \,\leq\,
\frac{d+1}{2}\,.
\end{equation}

\medskip \noindent
(3) Harris introduced a new idea to the irreducibility problem, which
completed the case of plane nodal curves (\cite{Harris}).  This new idea
was to proceed inductively from rational plane nodal curves (whose family is
classically known to be irreducible) to any family of plane nodal
curves of a given genus.  Further development of this idea lead to new
results by Ran \cite{Ran} and by Kang \cite{Kang2}: if $O_m$
denotes an ordinary singularity of order \mbox{$m\geq 2$}, then Ran showed
that, for each \mbox{$n\geq 0$}, the variety \mbox{$\Vdi(n\!\!\:\cdot\!\!\:A_1,
1\!\!\:\cdot \!\!\: O_m)$} is irreducible (or empty). Kang's result
says that, for each \mbox{$n\geq 0$}, \mbox{$k\leq 3$}, the variety \Vdnk\
is irreducible (if non-empty). 

However, the requirement to study all possible deformations of the considered
curves does not allow to extend such an approach to more complicated
singularities, or to a large number of singularities different from nodes.

\medskip \noindent
(4) Up to now, there is mainly one approach which is applicable to
equisingular families of curves of any degree with any quantity of arbitrary
singularities (and even to projective hypersurfaces of any dimension).
The basic idea is to find an irreducible analytic space $\km(S_1,\dots,S_r)$
and a dominant morphism
$$\VD \lra \km(S_1,\dots, S_r)$$
with equidimensional and irreducible fibres. It turns out, that in such a way
proving the irreducibility of $\VD$ can be reduced to an
$H^1$-vanishing problem:
let \mbox{$Z(C)=Z^s(C)$} (resp.\ \mbox{$Z(C)=Z^a(C)$}) be the
zero-dimensional schemes encoding the topological (resp.\ analytic)
type of the singularities (see Section \ref{sec:existence}).
Then the variety $\VD$ is irreducible if \mbox{$H^1\bigl(\kj_{Z(C)/\Sigma}(D)\bigr)=0$}
for each \mbox{$C\in \VD$}.

\medskip\noindent
For a detailed discussion of the latter
approach, we refer to \cite{GLS98,Kei}. Combining this approach with
\cite[Thm. 3]{GengXu95} and another $H^1$-vanishing theorem based on
the Castelnuovo function approach, we obtain:

\begin{theorem}[{\cite{GLS98}}]
\label{Theorem 3} Let \mbox{$S_1,\dots,S_r$} be topological or
analytic types of plane curve singularities,
and $d$ an integer.
If\/ \mbox{$\max\limits_{i=1..r} \tau'(S_i)\leq (2/5)d\!\!\:-\!\!\:1\,$} and
$$ \frac{25}{2}\cdot \# (\text{nodes\/})+18\cdot\#(
\text{cusps\/})+
\frac{10}{9} \cdot \!\!\!\sum_{\tau'(S_i)\ge 3}
\bigl(\tau'(S_i)\!\!\:+\!\!\:2\bigr)^2 \,<\: d^2,
$$
then $\Vd$ is non-empty and irreducible. Here, $\tau'(S_i)$ stands for
$\tes(S_i)$ if $S_i$ is a topological type and for $\tau(S_i)$ if
$S_i$ is an analytic type. 
\end{theorem}

\noindent
In particular,

\begin{corollary}
\label{nodes and cusps} Let \mbox{$d\geq 8$}. Then
\mbox{$\Vdi (n\!\!\:\cdot\!\!\: A_1,k\!\!\:\cdot\!\!\: A_2)$} is
irreducible if
\begin{equation}
\label{n & c}
\frac{25}{2}\!\; n + 18\!\: k \: < \: d^2.
\end{equation}
\end{corollary}

 \begin{corollary}
\label{ordinary singularities} Let $S_1,\dots,S_r$ be ordinary
singularities of order $m_1,\dots,m_r$, and assume that \mbox{$\max \,m_i \leq
(2/5)\!\; d$}. Then \mbox{$\Vdi (S_1,\dots,S_r)$} is non-empty and
irreducible if
\begin{equation}
\label{ord sin}
\frac{25}{2} \cdot \# (\text{nodes\/}) +
\sum_{m_i\geq 3} \frac{m_i^2(m_i\!\!\:+\!\!\:1)^2}{4} \: < \: d^2.
\end{equation}
\end{corollary}

\medskip\noindent
{\bf Reducible Equisingular Families.}
In \cite{GLS98}, we apparently gave the first series of reducible ESF
of plane cuspidal curves, where the different components cannot be
distinguished by the fundamental group of the complement of the
corresponding curves. If this happens, we say that the ESF has
components which are {\em anti-Zariski pairs}\index{anti-Zariski pair}. 

The following proposition gives infinitely many ESFs with anti-Zariski
pairs:

\begin{proposition}[{\cite{GLS98}}] \label{prop example}
Let $p,d$ be integers satisfying
\begin{equation}
p\ge 15,\quad
6p\!\;<\!\;d\!\;\le \!\;12p-\tfrac{3}{2}-\sqrt{35p^2\!-\!\!\:15p\!\!\;+
\!\!\;\tfrac{1}{4}}\ .
\label{eqnew}\end{equation}
Then the variety \mbox{$\Vdi(6p^2\!\cdot\!\!\: A_2)$} of irreducible plane
curves of degree $d$ with $6p^2$ cusps has components of different dimensions.

Moreover, for all curves \mbox{$\,C\in
  \Vdi(6p^2\!\cdot\!\!\: A_2)$} the fundamental group of the
complement is \mbox{$\pi_1(\P^2\!\!\:\setminus \!\!\;C)=\Z/d\Z$}.
\end{proposition}

\noindent
For instance, the variety \mbox{$V_{91}^{\text{\it irr}} (1350 \cdot A_2)$} is
reducible, and it has components which are anti-Zariski pairs.

Further, in \cite{GLS01}, we gave a series of reducible ESF of plane
curves with only ordinary singularities having components which are
anti-Zariski pairs:

\begin{proposition}[{\cite{GLS01}}]
\label{counterex}
Let \mbox{$m\geq 9$}. Then there is an integer
  \mbox{$\ell_0=\ell_0(m)$} such that for each \mbox{$\ell\geq
    \max\{\ell_0,m\}$} and for each $s$ satisfying
$$ \frac{\ell\!\!\;-\!\!\;1}{2} \:\leq  \:s\:\leq \:
\ell\left(1\!\!\:-\!\!\:\sqrt{\frac{2}{m}}\right)-\frac{3}{2} $$
the {variety\/} \mbox{$V_{\ell m+s}^{\text{\it irr}}(\ell^2\cdot O_m)$} of plane
irreducible curves of degree \mbox{$\ell m\!\!\:+\!\!\:s$} having
$\ell^2\!$ ordinary singularities of order $m$ as only 
singularities is reducible.

More precisely, \mbox{$V_{\ell m+s}^{\text{\it irr}}(\ell^2\!\cdot O_m)$} has at
least two components, one regular component (of the
expected dimension) and one component of higher dimension. And, for
each curve $C$ belonging to any of the components, we have
\mbox{$\pi_1(\P^2\setminus C)=\Z/(\ell m+s)\Z$}.
\end{proposition}

\section{Open Problems and Conjectures} \label{sec:problems}

Though some results discussed above are sharp, others seem to be far from
a final form, and here we start with a discussion and conjectures
about the expected progress in the geometry of families of singular curves.
Further discussion concerns possible generalizations of the methods
and open questions. 

\subsection{ESF of Curves}

\noindent
{\bf Existence of Curves with Prescribed Singularities}.
A natural question about the existence results for algebraic curves
given in Section \ref{sec:smoothness} concerns a possible improvement
of the asymptotically proper conditions to asymptotically optimal ones:

\medskip\noindent
{\it How to improve the
constant coefficients in the general sufficient conditions for the
existence\,?}

\smallskip\noindent Concerning our method based on $H^1$-vanishing
for the ideal sheaves, a
desired improvement would come from finding better $H^1$-vanishing
conditions for generic zero-dimensional schemes. For instance, from proving
the Harbourne-Hirschowitz conjecture, which gives (if true) the best possible
$H^1$-vanishing criterion for ideal sheaves of generic fat
schemes $Z(\bm)$.

\smallskip\noindent Another type of questions concerning curves
with specific singularities is the following:
the known sufficient and necessary conditions for the existence of
singular plane curves are formulated as bounds to sums of singularity
invariants. But it seems that there cannot be a general condition of
this type which is sufficient {\em and\/} necessary at the same
time. The simplest question of such kind is:

\medskip \noindent {\it Are there $k$ and $d$ such that a
curve of degree $d$ with $k$ cusps does exist, but with \mbox{$k'<k$} cusps
does not\,?}

\smallskip \noindent A candidate could be Hirano's series of
cuspidal curves mentioned in the introduction.

\bigskip \noindent
{\bf T-Smoothness and Versality of
Deformations}.\index{versality of deformation} The following conjecture about the asymptotic
properness of the sufficient conditions for the T-smoothness of
topological ESF of plane curves given in Section
\ref{sec:smoothness} seems to be quite realistic (and holds for
semiquasihomogeneous singularities):

\begin{conjecture}\label{cop2} There exists an absolute constant
  \mbox{$A>0$} such that for each
topological singularity type $S$ there are infinitely many pairs
\mbox{$(r,d)\in\N\!\:^2$}  such that
\mbox{$\Vdi(r\!\!\:\cdot\!\!\: S)$} is empty or non-smooth or has
dimension greater than the expected one and
\mbox{$r\cdot\gamma(S)\:\le\: A\cdot d^2$}.
\end{conjecture}

\noindent
We propose a similar conjecture for analytic ESF of plane curves, though
it is confirmed only for simple singularities (in which case it
coincides with the conjecture for topological ESF).

A closely related question, belonging to local singularity theory,
concerns the $\gam$-invariant:

\medskip\noindent
{\it Find an explicit formula, or an algorithm to compute $\ges(f)$,
$\gea(f)$. Find (asymptotically) close lower and upper bounds
for these invariants. Is $\ges$ a topological invariant\,?}

\bigskip \noindent {\bf Irreducibility Problem}. Our sufficient
irreducibility conditions seem to be far from optimal ones. We
state the problem:

\medskip \noindent
{\it Find asymptotically proper sufficient conditions for
the irreducibility of ESF of plane curves (or show that
the conditions in Section \ref{sec:irreducibility} are
asymptotically proper)}.

\medskip\noindent
We also rise the following important question:

\medskip\noindent
{\it Does there exist a pair of plane
irreducible algebraic curves of the same degree
with the same collection of singularities, which belong to
different components of an ESF but are topologically
isotopic in $\P^2$ (anti-Zariski pair)\,?}

\medskip\noindent The examples in Section \ref{sec:irreducibility}
provide candidates for this -- reducible ESF, whose members have
the same (Abelian) fundamental group of the complement.

\subsection{Hypersurfaces in Higher-Dimensional Varieties}

One
can formulate the existence, T-smoothness, and irreducibility
problems for families of hypersurfaces with isolated
singularities, belonging to (very) ample linear systems on
projective algebraic varieties. To find a relevant approach to
these problems is the most important question. For hypersurfaces of
dimension $>1$ there exists no infinitesimal deformation theory for
topological types. So, we restrict ourselves to analytic types here.

Constructions of curves with prescribed singularities as presented in Section
\ref{sec:existence} can, in principle, be generalized to higher
dimensions. An expected analogue of the results for curves could
be

\begin{conjecture}\label{cop3} Given a very ample linear system $|W|$
  on a projective algebraic
variety $X$ of dimension $n$, there exists a constant \mbox{$A=A(X,W)>0$}
such that, for each collection $S_1,\dots,S_r$ of singularity types
and for each positive integer $d$ satisfying
\mbox{$\sum_{i=1}^r\mu(S_i)<Ad^n$}, there is a hypersurface
\mbox{$W_d\in|dW|$} with exactly $r$ isolated singularities of
types $S_1,\dots,S_r$, respectively.
\end{conjecture}

\noindent In view of the patchworking approach, to prove the
conjecture, it is actually enough to consider the case of ordinary singularities
and to answer the following analogue of one of the above
questions affirmatively:

\medskip\noindent
{\it Does there exist some number \mbox{$A(n)>0$} such that, for each
analytic type $S$ of isolated hypersurface singularities in $\P^n$,
there exists a hypersurface of $\P^n$ of degree \mbox{$d\le 
A(n) \cdot \mu(S)^{1/n}$} which has a singularity of type $S$ and no
non-isolated singularities\,?}

\smallskip\noindent This is known only for simple singularities (see
\cite{West,Sh01}). 

\medskip\noindent
Hypersurfaces with specific singularities (such as nodes)
attract the attention of many researches, mainly looking for the
maximal possible number of singularities (see, for instance,
\cite{Chm}). We would like to raise the question about an analogue of the
Chiantini-Ciliberto theorem for nodal curves on surfaces (see \cite{CC})
as a natural counterpart, concerning the domain with regular behaviour
of ESF: 

\medskip\noindent
{\it
Given a projective algebraic variety $X$ and a very ample linear system
$|W|$ on it with a non-singular generic member. Prove that, for any
\mbox{$r\le\dim|W|$} there exists a hypersurface \mbox{$W_r\in|W|$}
with $r$ nodes as its only 
singularities such that the germ at $W$ of the corresponding ESF is
T-smooth.}

\subsection{Related Problems}

\noindent
{\bf Enumerative Problems}. Recently, the newly founded theories
of moduli spaces of stable curves and maps, Gromov-Witten
invariants, quantum cohomology, as well as deeply developed
methods of classical algebraic geometry and algebraic topology
have led to a remarkable progress in enumerative geometry,
notably for the enumeration of singular algebraic curves (see, for
example, \cite{KM,CH1,CH2,GP} for the enumeration of rational nodal
curves on rational surfaces, see \cite{Ran1,CH3} for the enumeration of
plane nodal curves of any genus, see \cite{Kaz,Liu} for counting
curves with arbitrary singularities). We point out that the
questions to which this survey has been devoted, such as on
the existence of certain singular algebraic curves, on the
expected dimension and on the transversality of the intersection of
ESF, are unavoidable in all of the above
approaches to enumerative geometry. The affirmative answers to
such questions are necessary for attributing an enumerative meaning
to the computations in the aforementioned works.

We pose the problem to find links between the methods discussed above
and the methods of enumerative geometry, and we expect that 
this would lead to a solution for new enumerative problems and to a better
understanding of known results. As an example, we mention the tropical
enumerative geometry \cite{Mi1,Mi2,Sh13}, in which the
patchworking construction and, more generally, the deformation theory
play an important role.

\bigskip \noindent {\bf Non-Isolated Singularities}. None of the
problems discussed above is even well-stated for
non-reduced curves, or hypersurfaces with non-isolated
singularities. We simply mention this as a direction for
further study.

\section*{Acknowledgments}

Work on the results presented in this paper has been supported by the
DFG Schwerpunkt ``Globale Methoden in der komplexen Geometrie''
and by the Hermann Minkowski -- Minerva Center for Geometry at Tel
Aviv University.




\printindex
\end{document}